\documentclass[11pt]{amsart}
\usepackage[]{amsmath, amsthm, amsfonts, verbatim, txfonts}
\DeclareFontFamily{OT1}{rsfs}{}
\DeclareFontShape{OT1}{rsfs}{n}{it}{<-> rsfs10}{}
\DeclareMathAlphabet{\mathscr}{OT1}{rsfs}{n}{it}

\newcommand\C{{\mathbb C}}
\newcommand\R{{\mathbb R}}

\def\cD{{\mathscr D}}
\def\cE{{\mathscr E}}

\def\ocM{\overline{\mathscr M}}
\def\oM{\overline M}

\def\cI{{\mathscr I}}
\def\cR{{\mathscr R}}
\def\cG{{\mathscr G}}
\def\oG{{\overline{\mathscr G}}}

\def\cT{{\mathscr T}}
\def\bR{{\mathbb R}}
\def\bZ{{\mathbb Z}}
\def\bQ{{\mathbb Q}}

\def\cT{{\mathscr T}}

\def\cZ{{\mathscr Z}}

\def\cP{{\mathscr P}}
\def\cC{{\mathcal C}}
\begin{document}
\newtheorem {theo}{Theorem}
\newtheorem {coro}{Corollary}
\newtheorem {lemm}{Lemma}
\newtheorem {rem}{Remark}
\newtheorem {defi}{Definition}
\newtheorem {ques}{Question}
\newtheorem {prop}{Proposition}
\def\spb{\smallpagebreak}
\def\mpb{\vskip 0.5truecm}
\def\bpb{\vskip 1truecm}
\def\wtM{\widetilde M}
\def\tM{\widetilde M}
\def\wtN{\widetilde N}
\def\tN{\widetilde N}
\def\tC{\widetilde C}
\def\tX{\widetilde X}
\def\tY{\widetilde Y}
\def\ti{\widetilde \iota}
\def\bs{\bigskip}
\def\ms{\medskip}
\def\ni{\noindent}
\def\td{\nabla}
\def\pd{\partial}
\def\hol{$\text{hol}\,$}
\def\Log{\mbox{Log}}
\def\bfQ{{\bf Q}}
\def\Todd{\mbox{Todd}}
\def\bP{{\bf P}}
\def\dxi{d x^i}
\def\dxj{d x^j}
\def\dyi{d y^i}
\def\dyj{d y^j}
\def\dzi{d z^I}
\def\dzj{d z^J}
\def\ozi{d{\overline z}^I}
\def\ozj{d{\overline z}^J}
\def\oz1{d{\overline z}^1}
\def\oz2{d{\overline z}^2}
\def\oz3{d{\overline z}^3}
\def\sI{\sqrt{-1}}
\def\hol{$\text{hol}\,$}
\def\ok{\overline k}
\def\ol{\overline l}
\def\oJ{\overline J}
\def\oT{\overline T}
\def\oS{\overline S}
\def\oV{\overline V}
\def\oW{\overline W}
\def\oI{\overline I}
\def\oK{\overline K}
\def\oL{\overline L}
\def\oj{\overline j}
\def\oi{\overline i}
\def\ok{\overline k}
\def\oz{\overline z}
\def\om{\overline mu}
\def\on{\overline nu}
\def\oa{\overline \alpha}
\def\ob{\overline \beta}
\def\of{\overline f}
\def\og{\overline \gamma}
\def\ogamma{\overline \gamma}
\def\odelta{\overline \delta}
\def\otheta{\overline \theta}
\def\ophi{\overline \phi}
\def\opd{\overline \partial}
\def\oA{\overline A} 
\def\oB{\overline B}
\def\oC{\overline C}
\def\op{\overline D}
\def\oIq1{\oI_1\cdots\oI_{q-1}}
\def\oIq2{\oI_1\cdots\oI_{q-2}}
\def\op{\overline \partial}
\def\ua{{\underline {a}}}
\def\us{{\underline {\sigma}}}
\def\tor{{\mbox{tor}}}
\def\vol{{\mbox{vol}}}
\def\rank{{\mbox{rank}}}
\def\Aut{{\mbox{Aut}}}
\def\Int{{\mbox{Int}}}
\def\Hol{{\mbox{Hol}}}
\def\Ad{{\mbox{Ad}}}
\def\dim{{\mbox{dim}}}
\def\bp{{\bf p}}
\def\bk{{\bf k}}
\def\a{{\alpha}}
\def\tchi{\widetilde{\chi}}
\def\b{{\bullet}}
\def\fg{{\mathfrak g}}
\def\ll{\lfloor}
\def\lr{\rfloor}
\title[Nonexistence of arithmetic fake compact hermitian symmetric spaces] 
{\sc Nonexistence of arithmetic fake compact hermitian symmetric spaces of type other than $A_n\: (n\leqslant 4)$}

\thanks{\ni \\
G. Prasad: University of Michigan, Ann Arbor, MI 48109\\
e-mail: gprasad@umich.edu\\
S.-K. Yeung:  Purdue University, West Lafayette, IN 47907\\
(corresponding author)
email: yeung@math.purdue.edu}
\maketitle
{\centerline{\sc Gopal Prasad and Sai-Kee Yeung}}

\vskip4mm

\noindent
{\bf Abstract.} {\it The quotient of a hermitian symmetric space of non-compact type by a torsion-free cocompact arithmetic subgroup of the identity component of the group of isometries of the symmetric space is called an 
arithmetic fake compact hermitian symmetric space
if it has the same Betti numbers as the compact dual of the hermitian symmetric space.  This is a natural generalization of the notion of ``fake projective planes ''
to higher dimensions.   Study of arithmetic fake compact hermitian symmetric spaces of type $A_n$ with even $n$ has been completed in
$\mathrm{[PY1]}$, $\mathrm{[PY2]}$.  The results of this paper, combined with those of $\mathrm [PY2]$, imply that there does not exist any arithmetic fake compact hermitian symmetric space of type other than  $A_n$, $n\leqslant 4$ (see Theorems 1 and 2 in the Introduction below and Theorem 2 of $\mathrm [PY2]$). The proof involves the volume formula given in $\mathrm{[P]}$, the Bruhat-Tits theory of reductive $p$-adic groups, and delicate estimates of various number theoretic invariants.}

\vskip3mm
\ni {\it Keywords:} arithmetic lattices, Bruhat-Tits theory, volume formula, cohomology.

\vskip3mm\ni {\it AMS 2010 Mathematics subject classification:} Primary 11F06, 22E40; Secondary 11F75 

\vskip6mm

\begin{center}
{\bf 1. Introduction}  
\end{center}
\vskip4mm

\ni{\bf 1.1.} Let $\oG$ be a connected real semi-simple Lie group with trivial center and with no nontrivial compact  normal subgroups, and $\fg$ be its Lie algebra. The group $\Aut(\oG)$ (=$\Aut(\fg)$) of automorphisms of $\oG$ is a Lie group with finitely many connected components, and  its identity component  is $\oG$. We will denote the identity component of $\Aut(\oG )$ in the {\it Zariski-topology}  by $\Int(\oG )$.  Let $X$ be the symmetric space of $\oG$ ($X$ is the space of maximal compact subgroups of $\oG$), and $X_u$ be the compact dual of $X$. There is a natural identification of the  group of isometries of $X$ with $\Aut(\oG )$. We  assume in this paper that $X$ (and hence $X_u$) is hermitian.  Then every holomorphic automorphism of $X$ is an isometry. The group $\Hol(X)$ of holomorphic automorphisms of $X$ is a subgroup of finite index of the group $\Aut(\oG)$ of isometries, and it is known (see [Ta], the remark in \S5)  that $\Hol(X)\cap \Int(\oG)=\oG$.  
\vskip1mm

\ni{\bf 1.2.} We will say that the quotient  $X/\Pi$ of $X$ by a {\it torsion-free} cocompact discrete subgroup $\Pi$ of $\oG$ is a {\it fake compact hermitian symmetric space}, or a {\it fake $X_u$}, if its Betti numbers are same as that of $X_u$; $X/\Pi$ is an {\it arithmetic fake compact hermitian symmetric space}, or an {\it arithmetic fake $X_u$}, if, moreover, $\Pi$ is  irreducible (i.e., no subgroup of $\Pi$ of finite index is a direct product of two infinite normal subgroups) and it is an arithmetic subgroup of $\oG$.  Any such space can be endowed with the structure of a smooth complex projective variety. \vskip1mm

We note that if $\oG$ contains an irreducible arithmetic subgroup, then the simple factors of its complexification are isomorphic to each other, see [Marg], Corollary 4.5 in Ch.\:IX. Also, if $\R$-rank of $\oG$ is at least $2$, which is the case for all $\oG$ to be considered in \S\S 4--7 of this paper, then by Margulis' arithmeticity theorem ([Marg], Ch.\:IX]), any irreducible discrete cocompact subgroup 
 of ${\oG}$ (in fact, any irreducible lattice) is arithmetic.

 \vskip1mm
 
 If $\Pi$ is a torsion-free cocompact discrete subgroup of $\oG$, then there is a natural embedding of $H^*(X_u, \C)$ in $H^*(X/\Pi, \C)$, see [B], 3.1 and 10.2, and hence $X/\Pi$ is a fake $X_u$ if and only if this embedding  is an isomorphism.

\vskip1mm

\ni{\bf 1.3.} Let $\oG$, $X$ and $X_u$ be as above, and let $\Pi$ be a torsion-free cocompact discrete subgroup of $\oG$. Let $Z = X/\Pi$.  If $Z$ is a fake $X_u$, then the Euler-Poincar\'e characteristic $\chi(Z)$ of $Z$, and so the Euler-Poincar\'e characteristic $\chi(\Pi)$ of  $\Pi$,  equals $\chi(X_u)$.  As $X$ has been assumed to be hermitian, the Euler-Poincar\'e characteristic of $X_u$ is positive. On the other hand, it follows from Hirzebruch proportionality principle, see [Ser], Proposition 23, that the Euler-Poincar\'e characteristic of $X/\Pi$ is positive if and only if the complex dimension of $X$ is even. Using the results of [BP], we can easily conclude that there are only finitely many irreducible arithmetic fake compact hermitian symmetric spaces of types other than $A_1$. It is of interest to determine them all.

\vskip1mm

\ni{\bf 1.4.} Hermitian symmetric spaces have been classified by \'Elie Cartan; see [H], Ch.\:IX. We recall that the noncompact  irreducible hermitian symmetric spaces are the symmetric spaces of Lie groups ${\rm SU}(n+1-m, m)$, 
 ${\rm SO}(2,2n-1)$, ${\rm Sp}({2n})$, ${\rm SO}(2, 2n-2)$, ${\rm SO}^*(2n)$, an absolutely simple real Lie group of type $E_6$ with Tits index $^2E^{16'}_{6,2}$, and  an absolutely simple real Lie group of type $E_7$ with Tits index $E^{28}_{7,3}$  (for Tits indices see Table II in [Ti1]).   The complex dimensions of these spaces are $(n+1-m)m$, $2n-1$, $n(n+1)/2$, $2n-2$, $n(n-1)/2$, $16$ and $27$ respectively.  The Lie groups listed above are of type $A_n$, $B_n$, $C_n$, $D_n$, $D_n$, $E_6$ and $E_7$ respectively.  We will say that a symmetric space is one of these types if it is a product of symmetric spaces of noncompact  simple Lie groups of that type, and say that a hermitian locally symmetric space is of one of these types if its simply connected cover is a hermitian symmetric space of that type.
\vskip1mm

The purpose of this paper is to prove the following two theorems.
\vskip2mm

\ni{\bf Theorem 1.} {\it There does not exist an irreducible arithmetic fake compact hermitian symmetric space of type other than $A_n$.} 
\vskip2mm

Regarding  spaces of type $A_n$, we have the following result. 
\vskip3mm

\ni{\bf Theorem 2.} {\it There does not exist an irreducible arithmetic fake compact hermitian symmetric space of type $A_n$ with $n>4$.}
\vskip1mm

\vskip1mm

%The following theorem is an immediate consequence of the above two theorems.
%\vskip2mm
%\noindent{\bf Theorem 3.} {\it Let $\oG$, $X$ and $X_u$ be as above and $\Pi$ be an irreducible cocompact torsion-free arithmetic subgroup of $\oG$. We assume that $X$ is hermitian.
%Then unless $X$ is of type $A_n$, with $n\leqslant 4$, the natural injective homomorphism $H^*(X_u,\C)\rightarrow H^*(X/\Pi,\C) = H^*(\Pi,\C)$ cannot be surjective.}   
% \vskip1mm

The proof of Theorem 1 is carried out in \S\S 4--7. Arithmetic fake compact hermitian symmetric spaces  of type $A_n$, with $n$ even, have been studied in detail in [PY1] and [PY2].  In [PY1] we have given a classification of ``fake projective planes", the first of which was constructed by David Mumford in [Mu] using $p$-adic uniformization. Note that fake projective planes are arithmetic fake compact hermitian symmetric spaces of type $A_2$. Using ingenious computer-assisted group theoretic computations, Cartwright and Steger ([CS])  have shown that  the twenty eight classes of fake projective planes of [PY1] altogether contain {\it fifty} distinct fake projective planes up to isometry with respect to the Poincar\'e metric [CS].  Since each of them supports
two distinct complex structures [KK,\S5], there are exactly {\it one hundred} fake  projective planes counted 
up to biholomorphism. In [PY2] we have shown that arithmetic fake compact hermitian symmetric spaces of type $A_n$, with $n$ even, can exist only for $n = 2, \,4,$ and have constructed four arithmetic fake ${\bf P}^4_{\C}$, four arithmetic fake Grassmannians  
${\bf Gr}_{2,5}$,  and five (irreducible) arithmetic fake ${\bf{P}}^2_{\C}\times{{\bf P}}^2_{\C}$. (Fake ${\bf P}^4_{\C}$ and fake ${\bf Gr}_{2,5}$ are of type $A_4$ and every fake ${\bf{P}}^2_{\C}\times{{\bf P}}^2_{\C}$ is of type $A_2$.) To prove Theorem 2 we therefore assume that $n$ is odd and $>3$. The proof occupies \S\S 8--9. We also prove some results for $n=3$, see Proposition 3 at the end of \S 8,  and 9.3.
\vskip1mm

In the following subsection we will explain the strategy of the proof, and fix notation which will be used throughout the paper. 

\vskip1mm

\ni{\bf 1.5.} Let $\oG$, $X$, $X_u$ be as in 1.1; $X $ will be assumed to be a hermitian symmetric space of one of the following types: $A_n$ with $n>3$ odd, $B_n$, $C_n$, $D_n$, $E_6$ and $E_7$. Assume, if possible, that $\oG$ contains a cocompact irreducible arithmetic subgroup $\Pi$ whose orbifold Euler-Poincar\'e characteristic $\chi(\Pi)$ equals $\chi(X_u)$. Then there exist a number field $k$, a connected adjoint absolutely simple algebraic $k$-group $\overline{G}$ of same type as $X$, real places $v_1, \ldots, v_r$ of $k$ such that $\overline{G}(k_v)$ is compact for every real place $v$ different from $v_1, \ldots, v_r$, $\oG$ is isomorphic to $\prod_{j=1}^r\overline{G}(k_{v_j})^{\circ}$ (and will be identified with it),  and $\Pi$ is an arithmetic subgroup contained in $\overline{G}(k)$.  It is obvious from this that $k$ is totally real.  If $\overline{G}$ is either of type $A_n$ ($n>1$ arbitrary) or $D_n$ with $n$ odd, or of type $E_6$, then for every real place $v$ 
of $k$, $\overline{G}$ is an outer form over $k_v$, and hence the unique quadratic extension $\ell$ of $k$ over which $\overline{G}$ is an inner form is totally complex. If $\overline{G}$ is of type $D_n$ with $n$ even, and it is {\it not} a triality form of type $D_4$, then at every real place $v$ of $k$, $\overline{G}$ is an inner form, and hence either $\overline{G}$ is an inner $k$-form, or the unique quadratic extension $\ell$ over which $\overline{G}$ is an inner form is a totally real field. If $\overline{G}$ is a triality form of type $D_4$, let $\ell$ be a fixed cubic extension of $k$ contained in the smallest Galois extension of $k$ over which $\overline{G}$ is an inner form. For triality forms occuring in this paper, $\ell$ is totally real. As $\Pi$ is cocompact, by Godement compactness criterion $\overline{G}$ is anisotropic over $k$ (i.e., its $k$-rank is $0$).
\vskip1mm

Let $\pi: G \rightarrow \overline{G}$ be the simply connected cover of $\overline{G}$ defined over $k$. The kernel of the isogeny $\pi$ is the center $C$ of the simply connected $k$-group $G$.

\vskip1mm

{\it Description of $C$:} For a positive integer $a$, let $\mu_a$ be the kernel of the endomorphism $x\mapsto
x^a$ 
of ${\rm  GL}_1$. Then if $G$ is of type $^2A_n$, its center is $k$-isomorphic to the kernel of the norm map $N_{\ell/k}:\,R_{\ell/k}(\mu_{n+1}) \to \mu_{n+1}$, and if $G$ is of type $^2E_6$, its center is $k$-isomorphic to the kernel of the norm map 
$N_{\ell/k}:\,R_{\ell/k}(\mu_3) \to \mu_3$. If $G$ is of type $B_n$, $C_n$ or $E_7$, then $C$ is $k$-isomorphic to $\mu_2$. If $G$ is an inner $k$-form of type $D_n$ with $n$ even, then $C$ is $k$-isomorphic to $\mu_2\times \mu_2$, and if $G$ is a non-triality outer form of type $D_n$, $C$ is $k$-isomorphic to $R_{\ell/k}(\mu_2)$ or to the kernel of the norm map $N_{\ell/k}:\,R_{\ell/k}(\mu_4)\to \mu_4$ according as $n$ is even or odd.  If $G$ is a triality form of type $D_4$, let the cubic extension $\ell$ of $k$ be as above. Then $C$ is $k$-isomorphic to the kernel of the norm map $N_{\ell/k}:\,R_{\ell/k}(\mu_2)\to \mu_2$. 
\vskip1mm

It is known, and easy to see using the above description of $C$, that for any real place $v$ of $k$, the order of the kernel of the induced homomorphism $G(k_v)\rightarrow \overline{G}(k_v)$ is  $n+1$ if $G$ is of type $^2A_n$, is of order $2$ if $G$ is of type $B_n$, $C_n$ or $E_7$, is of order $4$ if $G$ is of type $D_n$, and of order $3$ if it is of type $^2E_6$. Moreover, as $G(k_v)$ is connected, $\pi(G(k_v)) = \overline{G}(k_v)^{\circ}$. 
Let  $\cG =\prod_{j=1}^rG(k_{v_j})$, and let $\widetilde\Pi$ be the inverse image of $\Pi$ in $\cG$. Then the kernel of the homomorphism $\pi : \cG\to \oG$ is of order $s^r$, and hence the orbifold Euler-Poincar\'e characteristic $\chi(\widetilde\Pi)$ of $\widetilde\Pi$ equals $\chi(\Pi)/s^r = \chi(X_u)/s^r$, where here, and in the sequel, $s = n+1$ if $G$ is of type $A_n$, $s =2$ if $G$ is of type $B_n$, $C_n$ or $E_7$, $s =4$ if $G$ is of type $D_n$, and $s = 3$ if $G$ is of type $E_6$. Now let $\Gamma$ be a maximal discrete subgroup of $\cG$ containing $\widetilde\Pi$. Then the orbifold Euler-Poincar\'e characteristic $\chi(\Gamma)$ of $\Gamma$ is a submultiple{\footnote {given two nonzero real numbers $x$ and $y$, we say that $y$ is a {\it submultiple} of $x$ if $x/y$ is an integer}} of $\chi(\widetilde\Pi)=\chi(X_u)/s^r$.  Using the volume formula of [P], some nontrivial number theoretic estimates, the Bruhat-Tits theory, and the Hasse principle for semi-simple groups (Proposition 7.1 of  [PR]), we will show  that $\cG$ does not contain such a subgroup $\Gamma$. This will prove Theorems 1 and 2.

\vskip4mm

\centerline{\bf 2. Preliminaries}

\vskip3mm

\ni{\bf 2.1.} We will use the notations introduced in 1.5. Thus $k$ will be a totally real number field, $G$ an absolutely simple simply connected algebraic $k$-group (of one of the following nine  types: $^2A_n$ with $n\,(>3)$ odd, $B_n$, $C_n$, $^1D_n$, $^2D_n$, $^3D_4$, $^6D_4$, $^2E_6$, and $E_7$), $C$ its center, $\cG = \prod_{j=1}^rG(k_{v_j})$. We will think of $G(k)$ as a subgroup of $\cG$ in terms of its diagonal embedding. 
\vskip1mm

 $V_f$ (resp.\,\,$V_{\infty}$) will denote the set of nonarchimedean
(resp.,\,archimedean) places of $k$. As $k$ admits at least $r$ distinct real places, see 1.5, $d :=[k:\bQ]\geqslant r$.  For $v\in V_f$, $q_v$ will
denote 
the cardinality of 
the residue field ${\mathfrak f}_v$ of $k_v$.  If $G$ is an outer form, $\ell$ will denote the quadratic or cubic extension of $k$ as in 1.5. If $G$ is an inner form, let $\ell = k$. 

\vskip1mm
As explained in 1.5, to prove Theorems 1 and 2 it will suffice to show that $\cG$ does not contain a maximal arithmetic subgroup 
$\Gamma$ ($\Gamma$ arithmetic with respect to the $k$-structure on $G$) whose orbifold Euler-Poincar\'e characteristic is a submultiple of $\chi(X_u)/s^r$. Assume, if possible, that such a $\Gamma$ exists. Then $\Lambda : = \Gamma\cap G(k)$ is a ``principal" arithmetic subgroup, i.e., for every nonarchimedean place $v$ of $k$, the closure $P_v$ of $\Lambda$ in $G(k_v)$ is a parahoric subgroup  and  $\Lambda = G(k)\cap \prod_{v\in V_f}P_v$, moreover, $\Gamma$ is the normalizer of $\Lambda$ in $\cG$; see Proposition 1.4(iv) of [BP]. Let the ``type" $\Theta_v$ of $P_v$ be as in [BP], 2.2, and $\Xi_{\Theta_v}$ be as in 2.8 there. If $P_v$ is hyperspecial, then $\Xi_{\Theta_v}$ is trivial. The order of $\Xi_{\Theta_v}$ is always a divisor of $s$ ($s$ as in 1.5).  We note that for all but finitely many $v\in V_f$, $P_v$ is hyperspecial. 
\vskip1mm

In terms 
of the normalized Haar-measure $\mu$ on $\cG= \prod_{j=1}^r G(k_{v_j})$ used in [P]
and [BP], and to be used in this paper, $\vert\chi(\Gamma)\vert=\chi(X_u)\mu(\cG/\Gamma)$ (see [BP], 4.2). Thus the condition that $\chi(\Gamma)$ is a submultiple of  $\chi(X_u)/s^r$ is equivalent to the condition that $\mu(\cG/\Gamma)$ is a submultiple of ${1/s^r}.$ We will show below that $\cG$ does not contain a maximal arithmetic subgroup $\Gamma$ such that $\mu(\cG/\Gamma)$ is a submultiple of ${1/s^r}$.
\vskip2mm

For a comprehensive survey of the basic notions and the main results of the Bruhat-Tits theory of reductive groups over nonarchimedean local fields, used in this paper, see [Ti2].
\vskip2mm 

\ni{\bf 2.2.} All unexplained notations are as in [BP] and [P]. Thus for a number field $K$, $D_K$ will denote the absolute value of its 
discriminant, $h_K$ its class number, i.e., the order of its class group 
$Cl (K)$. We will denote by $h_{K,s}$ the order of the subgroup of $Cl(K)$ consisting 
of the elements of order dividing $s$, where, as in 1.5, $s = n+1$ if $G$ is of type $A_n$, $s =2$ if $G$ is of type $B_n$ or $C_n$, $s =4$ if $G$ is of type $D_n$, and $s=3$ if $G$ is of type $E_6$. Then $h_{K,s}
\vert h_K$. We will denote by $U_K$ the multiplicative-group of units of $K$, and
by $K_s$ the subgroup of $K^{\times}$ consisting of the elements $x$
such that for every normalized valuation $v$ of $K$, $v(x)\in s\bZ$.
\vskip1mm
 
 \vskip1mm

\ni{\bf 2.3.} For a parahoric subgroup $P_v$ 
of $G(k_v)$,  we define $e(P_v)$ and $e'(P_v)$ by the following formulae (cf.\,Theorem 3.7 of [P]):
\begin{equation}
e(P_v)=\frac{q_v^{(\dim\,\oM_v+\dim\,\ocM_v)/2}}{\#\oM_v({\mathfrak f}_v)}.
\end{equation} 
\begin{equation}
e'(P_v) =e(P_v)\cdot \frac{\#\ocM_v({\mathfrak f}_v)}{q_v^{\dim\,\ocM_v}}=q_v^{(\dim\,\oM_v -\dim\,\ocM_v)/2}\cdot\frac{\#\ocM_v({\mathfrak f}_v)}{\#\oM_v({\mathfrak f}_v)}.
\end{equation}
\vskip1mm

\ni{\bf 2.4.} Let $m_1,\dots , m_n$ ($m_1\leqslant \cdots \leqslant  m_n$), where $n$ is the absolute rank of $G$, be the exponents of the Weyl group of $G$. For type $A_n$, $m_j = j$; for types $B_n$ and $C_n$, $m_j = 2j-1$;  for type  $D_n$ the exponents are $1$, $3$, $5$, \dots , $2n-5$, $2n-3$ and $n-1$ (the multiplicity of $n-1$ is two when $n$ is even);  for type $E_6$, the exponents are $1$, $4$, $5$, $7$, $8$ and $11$; and for type $E_7$, the exponents are $1$, $5$, $7$, $9,$  $11$, $13$ and $17$.  Then
\vskip2mm

\ni$\bullet$ if either $G$ is of inner type, or $v$ completely splits in $\ell$,  
$$e'(P_v) = e(P_v)\prod_{j=1}^{n}\Big(1-\frac{1}{q^{m_j+1}_v}\Big);$$ 

\vskip1mm

\ni$\bullet$ if $v$ does not split in $\ell$ and  $G$ is of type $^2A_n$ with $n$ odd, then 
$$e'(P_v) = e(P_v)\Big(1-\frac{1}{q_v^{n+1}}\Big)\prod_{j=1}^{(n-1)/2}\Big(1-\frac{1}{q_v^{2j}}\Big)\Big(1+\frac{1}{q_v^{2j+1}}\Big),$$ or 
$$e'(P_v) =e(P_v)\prod_{j=1}^{(n+1)/2}\Big(1-\frac{1}{q_v^{2j}}\Big)$$ according as $v$ {\it does not} or {\it does} ramify in $\ell$. 
\vskip1mm

\ni$\bullet$ if $G$ is of type $^2D_n$ and  $v$ does not split in $\ell$,    
$$e'(P_v) = e(P_v)\Big(1+\frac{1}{q_v^n}\Big)\prod_{j=1}^{n-1}\Big(1-\frac{1}{q_v^{2j}}\Big),$$ or $$e'(P_v) = e(P_v)\prod_{j=1}^{n-1}\Big(1-\frac{1}{q_v^{2j}}\Big)$$
according as $v$ {\it does not} or {\it does} ramify in $\ell$. 

\vskip1mm

\ni$\bullet$ if $G$ is a triality form (i.e., of type $^3D_4$ or $^6D_4$) and $v$ does not completely split in $\ell$, let $\omega$ be a nontrivial cube root of unity, then 
\vskip1mm

{(i)} if $\ell _v := \ell\otimes_k k_v$ is a (cubic) field extension of $k_v$, $$e'(P_v) = e(P_v)\Big(1-\frac{1}{q_v^2}\Big)\Big(1-\frac{\omega}{q_v^4}\Big)\Big(1-\frac{\omega^2}{q_v^4}\Big)\Big(1-\frac{1}{q_v^6}\Big),$$  or $$e'(P_v) = e(P_v) \Big(1-\frac{1}{q_v^2}\Big)\Big(1-\frac{1}{q_v^6}\Big)$$ according as $\ell_v$ is a {\it unramified} or a {\it ramified} extension of $k_v$,
\vskip1mm

(ii) if $\ell\otimes_k k_v$ is a direct product of $k_v$ and a quadratic field extension of $k_v$, then $$e'(P_v) = e(P_v) \Big(1+\frac{1}{q_v^4}\Big)\Big(1-\frac{1}{q_v^2}\Big)\Big(1-\frac{1}{q_v^4}\Big)\Big(1-\frac{1}{q_v^6}\Big),$$ or $$ e'(P_v) = e(P_v) \Big(1-\frac{1}{q_v^2}\Big)\Big(1-\frac{1}{q_v^4}\Big)\Big(1-\frac{1}{q_v^6}\Big)$$ according as the quadratic extension is {\it unramified} or {\it ramified}
\vskip1mm

\ni$\bullet$ if $G$ is of type $^2E_6$ and  $v$ does not split in $\ell$,    
$$e'(P_v) = e(P_v)\Big(1-\frac{1}{q_v^2}\Big)\Big(1+\frac{1}{q_v^5}\Big)\Big(1-\frac{1}
{q_v^6}\Big)\Big(1-\frac{1}{q_v^8}\Big)\Big(1+\frac{1}{q_v^9}\Big)\Big(1-\frac{1}{q_v^{12}}\Big),$$ or $$e'(P_v) = e(P_v)\Big(1-\frac{1}{q_v^2}\Big)\Big(1-\frac{1}
{q_v^6}\Big)\Big(1-\frac{1}{q_v^8}\Big)\Big(1-\frac{1}{q_v^{12}}\Big)
$$
according as $v$ {\it does not} or {\it does} ramify in $\ell$. 

\vskip1mm
 
\ni{\bf 2.5.} Since $q_v^{\dim\,\ocM_v}> \#\ocM_v({\mathfrak f}_v)$ (cf.\,2.6 of [P]), $e'(P_v) < e(P_v)$. It is not difficult to check by case-by-case computations, using (2) and the Bruhat-Tits theory,  that {\it for all $v\in V_f$, and an arbitrary  parahoric subgroup $P_v$ of $G(k_v)$,  $e'(P_v)$ is an integer.}  If, for example, either $G$ is quasi-split over $k_v$ and splits over the maximal unramified extension of $k_v$ (equivalently, $G(k_v)$ contains a hyperspecial  parahoric subgroup), or it does not split over the maximal unramified extension of $k_v$, then explicit computations can be avoided using the fact that the order of a subgroup of a finite group divides the order of the latter, an analogue (see [Gi]) for reductive groups over finite fields of a result of Borel and de Siebenthal on subgroups of maximal rank of a compact Lie group, and the fact that over a finite field $\mathfrak{f}$, the groups of $\mathfrak{f}$-rational points of connected absolutely simple $\mathfrak{f}$-groups of types $B_m$ and $C_m$, for an arbitrary $m$, have equal order. A detailed proof of the integrality of $e'(P_v)$ for groups of type $A_n$ is given in [GM].     
\vskip1mm

\ni{\bf 2.6.} Now we will use the volume formula of [P] to write down
the precise value of $\mu(\cG/\Lambda)$. As the Tamagawa number 
$\tau_k(G)$ of $G$ equals $1,$ Theorem 3.7 of [P] (recalled in 3.7 of [BP]), for $S=V_{\infty}$,  provides us the following:
\begin{equation}
\mu(\cG /\Lambda )= D_k^{\frac{1}{2}{\dim\,G}}(D_{\ell}/D_k^{[\ell:k]})^{\frac{1}{2}{\mathfrak s}}\Big(\prod_{j=1}^{n}\frac{m_j!}{(2\pi)^{m_j+1}}\Big)^{d}\cE, 
\end{equation}
where
$n$ is the absolute rank of $G$, $\mathfrak s$ is $(n-1)(n+2)/2$ if $G$ is of type $^2A_n$ with $n$ odd,  $2n-1$ if $G$ is of type $^2D_n$, $7$ if $G$ is a triality form (i.e., of type $^3D_4$ or $^6D_4$), $26$ if $G$ is of type $^2E_6$, and $0$ for all other groups under consideration in this paper, and $$\cE=\prod_{v\in V_f}e(P_v),$$ with $e(P_v)$ as in 2.3.
 \vskip1mm
  
\ni{\bf 2.7.} Let $\zeta_k$, $\zeta_{\ell}$ be the Dedekind zeta-functions of $k$ and $\ell$ respectively. We will let $\zeta_{\ell |k}$ denote the function $\zeta_{\ell}/\zeta_k$. If $\ell$ is a quadratic extension of $k$, which will often be the case in this paper, $\zeta_{\ell |k}$ is the Hecke $L$-function associated to the nontrivial Dirichlet character of $\ell/k$. Recall that 
$$\zeta_k(a) =\prod_{v\in V_f}\Big(1-\frac1{q_v^a}\Big)^{-1},$$ and if  $\ell$ is a quadratic extension of $k$, 
$$\zeta_{\ell|k}(a) ={\prod} '\Big(1-\frac1{q_v^a}\Big)^{-1}{\prod} ''
\Big(1+\frac1{q_v^a}\Big)^{-1},$$
where $\prod'$ is the product over the nonarchimedean places $v$ of $k$ which 
split in $\ell$, and $\prod''$
is the product over the nonarchimedean places $v$ which do not 
split and also do not ramify in $\ell$. We will let the reader write down a similar product expression for $\zeta_{\ell | k}(a)=\zeta_{\ell}(a)/\zeta_k(a)$ when $\ell$ is a cubic extension of $k$. 
\vskip1mm

Using the values of $e'(P_v)$ given in 2.4 we will rewrite the Euler product $\cE$ appearing in (3). For this purpose we define 
$$\cZ= \prod_{j=1}^{n} 
\zeta_k(m_j+1)$$
if $G$ is of inner type;
$$\cZ = \zeta_k(n+1)\prod_{j=1}^{(n-1)/2}\zeta_k(2j)\zeta_{\ell |k}(2j+1)$$
if $G$ is of type $^2A_n$ with $n$ odd; 
$$\cZ= \zeta_{\ell | k}(n)\prod_{j=1}^{n-1} 
\zeta_k(2j)$$ if $G$ is of type $^2D_n$; $$\cZ = \zeta_k(2)\zeta_{\ell |k}(4)\zeta_k(6)$$ if $G$ is a triality form; $$\cZ =\zeta_k(2)\zeta_{\ell |k}(5)\zeta_k(6)\zeta_k(8)\zeta_{\ell |k}(9)\zeta_k(12)$$ if $G$ is of type $^2E_6$. Then for all $G$, \begin{equation}
\cE = \cZ\prod_{v\in V_f} e'(P_v).\end{equation}

\ni{\bf 2.8.} If $G$ is of inner type, let
 \begin{equation}
\cR = 2^{-dn}|\prod_{j=1}^{n}   \zeta_k(-m_j)|. 
\end{equation}

If $G$ is of type $^2A_n$ with $n$ odd, let 
\begin{equation}
\cR = 2^{-dn}|\zeta_k(-n)\prod_{j=1}^{(n-1)/2}\zeta_k(1-2j)\zeta_{\ell |k}(-2j)|.
\end{equation}

If $G$ is of type $^2D_n$  (recall that $\ell$ is totally real or totally complex according as $n$ is even or odd), let 
\begin{equation}
\cR = 2^{-dn}|\zeta_{\ell |k}(1-n)\prod_{j =1}^{n-1}\zeta_k(1-2j)|.
\end{equation}
If $G$ is a triality form (then $\ell$ is a totally real cubic extension of $k$), let
\begin{equation}
\cR =2^{-4d}|\zeta_k(-1)\zeta_{\ell |k}(-3) \zeta_k(-5)| .
\end{equation}
If $G$ is of type 
$^2E_6$ (then $\ell$ is a totally complex quadratic extension of $k$), let
\begin{equation}
\cR = 2^{-6d}|\zeta_k (-1)\zeta_{\ell | k}(-4)\zeta_k (-5)\zeta_k (-7)\zeta_{\ell | k}(-8)\zeta_k (-11)|.
\end{equation}
\vskip1mm

Using the following functional equations for any totally real $k$ and respectively a totally real extension of arbitrary degree and a totally complex quadratic extension $\ell$ of $k$, $$\zeta_k(2a) = D_k^{\frac{1}{2}-
2a}\Big(\frac{(-1)^{a}2^{2a-1}\pi^{2a}}{(2a-1)!}\Big)^d \zeta_k(1-2a),$$
$$\zeta_{\ell |k}(2a) =\Big(\frac{D_k}{D_{\ell}}\Big)^{2a-\frac{1}{2}}\Big(\frac{(-1)^a2^{2a-1}\pi^{2a}}{(2a-1)!}\Big)^{d([\ell :k]-1)}\zeta_{\ell | k}(1-2a),$$ and  $$\zeta_{\ell|k}(2a+1) = \Big(\frac{D_k}{D_{\ell}}\Big)^{2a+\frac{1}{2}}\Big(\frac{(-1)^a 2^{2a}\pi^{2a+1}}{(2a)!}\Big)^d \zeta_{\ell|k}(-2a),$$
for every positive integer $a$, and the fact that $\dim\,G = n+2\sum m_j$, the volume formula (3) and the explicit value of $\cE$ given for each case  in 2.7, we find that 
\begin{equation}
\mu(\cG/\Lambda) = \cR \prod_{v\in V_f}e'(P_v),
\end{equation}
where $\cR$ is as above.
\vskip2mm

\ni{\bf 2.9.}  We have the following\begin{equation}
\mu(\cG/\Gamma)= \frac{\mu(\cG/\Lambda)}{[\Gamma :\Lambda ]}= \frac{\cR \prod_{v\in V_f}e'(P_v)}{[\Gamma :\Lambda]}.
\end{equation}
Let $s$ be as in 1.5.  Proposition 2.9 of [BP] applied to $G' = G$ and $\Gamma' = \Gamma$ implies that any prime  divisor of $[\Gamma :\Lambda]$ divides $s$. Now since $e'(P_v)$ is an integer for all $v\in V_f$, we 
conclude from (11) that  if $\mu(\cG/\Gamma)$ is a submultiple of $1$, then any prime which divides the numerator of the rational number 
$\cR$ is a divisor of $s$. We record this observation as the following proposition.
\begin{prop} If $\mu(\cG/\Gamma)$ is a submultiple of $1\,($or, equivalently,  the orbifold Euler-Poincar\'e characteristic  
\,$\chi(\Gamma)$ of \,$\Gamma$  is a submultiple of 
$\chi(X_u))$, then every prime divisor of the numerator of the rational number $\cR$ divides $s$.
\end{prop}

 \ni{\bf 2.10.} Let $\cT$ be
the set of all nonarchimedean places $v$ of $k$ such that {\it either} (i) $v$ does not ramify in $\ell$ (equivalently, $G$ splits over the maximal unramified extension of $k_v$) and  
$P_v$ is not a hyperspecial parahoric subgroup of $G(k_v)$, {\it or} (ii) $v$ ramifies in $\ell$, $G$ is quasi-split over $k_v$ and  $P_v$ is not special.  It can be easily seen, using the relative local Dynkin diagram of $G/k_v$ given in 4.3 of [Ti2], that if $v\notin \cT$, then $\Xi_{\Theta_v}$ is trivial; if $v\in\cT$ ramifies in $\ell$, then $\#\Xi_{\Theta_v}\leqslant 2$.  
\vskip1mm

If for a $v\in V_f$, $P_v$ is hyperspecial, then obviously  $e'(P_v) = 1$. On the other hand, it is not difficult to see, by  direct computation, that $e'(P_v)>s$ for all $v\in \cT$.  Therefore,  
$\cE =\prod_{v\in V_f} e(P_v)>\prod_{v\in V_f} e'(P_v)>s^{\#\cT}$ (cf.\,2.5), and hence, we see from (3) that 
\begin{equation}
\mu(\cG/\Lambda )> D_k^{\frac{1}{2}{\dim\,G}}(D_{\ell}/D_k^{[\ell:k]})^{\frac{1}{2}{\mathfrak s}}\Big(\prod_{j=1}^{n}\frac{m_j!}{(2\pi)^{m_j+1}}\Big)^{d}s^{\#\cT}.
\end{equation}
Since $\mu(\cG/\Gamma) =\mu(\cG/\Lambda)/[\Gamma:\Lambda]$ is a submultiple of $1/s^r$ (see 2.1), we conclude that $\mu(\cG/\Lambda) \leqslant [\Gamma:\Lambda]/s^r$. From bound (12) we now obtain:
\begin{equation}
 D_k^{\frac{1}{2}{\dim\,G}}(D_{\ell}/D_k^{[\ell:k]})^{\frac{1}{2}{\mathfrak s}}\Big(\prod_{j=1}^{n}\frac{m_j!}{(2\pi)^{m_j+1}}\Big)^{d}s^{\#\cT}< [\Gamma:\Lambda]/s^r.
 \end{equation}

\bs

\centerline {\bf 3.  Discriminant bounds}
\bs

We will recall discriminant bounds required in later discussions.
We define $M_r(d)=\min_K D_K^{1/d},$ where the minimum is taken over all
totally real number fields $K$ of degree $d.$  Similarly, we define  $M_c(d)=\min_K D_K^{1/d},$ by taking the minimum over all
totally complex number fields $K$ of degree $d.$
 
The precise values of $M_r(d), M_c(d)$ for low values of $d$  are given in the following table (cf.\:[N]).  

$$\begin{array}{cccccccc}
d:&2&3&4&5&6&7&8\\
M_r(d)^d:
&5
&49
&725
&14641
&300125
&20134393
&282300416\\
M_c(d)^d:&3&&117&&9747&&1257728.
\end{array}
$$

\vskip1mm
The following proposition can be proved in the same way as Proposition 2  
in [PY2] has been proved.

\begin{prop}
Let $k$ and $\ell$ be a totally real number field and a
totally complex number field of degree $d$ respectively.    
$$\begin{array}{cccccccc}
\forall d\:\geqslant&2&3&4&5&6&7&8\\

D_k^{1/d}>&2.23&3.65&5.18&6.8&8.18&11.05&11.38\\
\ \ 

D_{\ell}^{1/d}>&1.73&&3.28&&4.62&&5.78.
\end{array}$$
\end{prop}

\bs
\begin{center}
{\bf 4. $G$ of type ${B}_n$ or ${C}_n$}  
\end{center}
\vskip4mm

\ni{\bf 4.1.}  In this section we assume that $G$ is of type $B_n$ or $C_n$ with $n>1$. Then its dimension is $n(2n+1)$. The $j$-th exponent $m_j = 2j-1$, $\mathfrak{s}=0$, and the complex dimension of the symmetric space $X$ of $\cG =\prod_{j=1}^r G(k_{v_j})$ is $r(2n-1)$ if $G$ is of type $B_n$, and is $rn(n+1)/2$ if $G$ is of type $C_n$. The center $C$ of $G$ is $k$-isomorphic to $\mu_2$ and $s =2$. 
The Galois cohomology group $H^1(k, C)$ is isomorphic to $k^{\times}/{k^{\times}}^{2}$. The order of the first term of the short exact sequence of Proposition 2.9 of [BP], for $G' =G$ and $S = V_{\infty}$, is $2^{r-1}$. From the proof of Proposition 0.12 of [BP], we easily conclude that $\#k_2/{k^{\times}}^2 \leqslant h_{k,2}2^d$. Let $\cT$ be as in 2.10. We can adapt the argument used to prove Proposition 5.1 in [BP], and the argument in 5.5 of [BP], for $S=V_{\infty}$ and $G' =G$, to derive the following bound from Proposition 2.9 of [BP]:
\begin{equation}
[\Gamma : \Lambda]\leqslant h_{k,2}2^{d+r-1+\#\cT}.
\end{equation}

\ni Hence, from (13) we obtain 
 \begin{equation}
D_k^{1/d} <
f_1(n,d,h_{k,2}):=\Big[ \{
2\prod_{j=1}^{n}\frac{(2\pi)^{2j}}{(2j-1)!}\}^d\cdot \frac{h_{k,2}}2
\Big]^{\frac{2}{dn(2n+1)}}.
\end{equation}

According to the Brauer-Siegel Theorem, for a totally real number field $k$ of degree $d$, and all real $\delta >0$, 
\begin{equation}
h_k R_k\leqslant 2^{1-d}{\delta (1+\delta )}\Gamma((1+\delta)/2)^d(\pi^{-d}D_k)^{\frac{(1+\delta)}{2}}
\zeta_k(1+\delta),
\end{equation}
where $R_k$ is the regulator of $k$. Now from (15) we get the following bound:
\begin{equation}
{D_k^{1/d}<f_2(n,d,R_{k},\delta)\hskip8cm}
\end{equation}

$${\hskip1.3cm}:=\Big[\{\frac{\Gamma((1+\delta)/2)\zeta(1+\delta)}{\pi^{(1+\delta)/2}}
\prod_{j=1}^{n}\frac{(2\pi)^{2j}}{(2j-1)!}
\}\cdot
\{\frac{\delta(1+\delta)}{R_{k}}\}^{\frac{1}{d}} \Big]^{\frac{2}{(2n^2+n-1-\delta)}},$$
since $\zeta_k(1+\delta)\leqslant \zeta(1+\delta)^d$, where $\zeta = \zeta_{\bQ}$. Using the lower bound  
$R_k\geqslant 0.04\,e^{0.46d}$, for a totally real number field $k$, due to R.\,Zimmert [Z], we obtain the following bound from (17):
\begin{equation}
{D_k^{1/d} <f_3(n,d,\delta) \hskip8cm}
\end{equation}
$${\hskip2cm :=\Big[\{\frac{\Gamma((1+\delta)/2)\zeta(1+\delta)}{\pi^{(1+\delta)/2}e^{0.46}}
\prod_{j=1}^{n}\frac{(2\pi)^{2j}}{(2j-1)!}
\}\cdot
\{25\delta(1+\delta)\}^{\frac{1}{d}} \Big]^{\frac{2}{(2n^2+n-1-\delta)}} .}$$

\bs
\ni{\bf 4.2.}  It is obvious that for fixed $n\geqslant 2$ and $\delta \in [0.04,\, 9]$, $f_3(n,d,\delta)$ decreases as $d$ increases.  Now we observe that for $n\geqslant 9$, $(2n-1)!> (2\pi)^{2n}$. From this it is easy to see that if for a given $d$, $\delta\in [0.04, 9]$, and $n \geqslant 8$, $f_3(n,d,\delta) \geqslant 1$, then $f_3(n+1, d,\delta)< f_3(n,d,\delta)$, and if $f_3(n,d,\delta)< 1$, then $f_3(n+1,d,\delta)<1$.  In particular, if for given $d$, and $\delta\in [0.04,\, 9]$, $f_3(8, d,\delta)<c$, with $c\geqslant 1$, then $f_3(n,d',\delta)<c$ for all $n\geqslant 8$ and $d'\geqslant d$.  
\vskip1mm

We obtain by a direct computation 
the following upper bound for
the value of $f_3(n,2,3)$ for $6\leqslant n\leqslant 14$.

$$\begin{array}{cccccccccc}
n:&14&13&12&11&10&9&8&7&6\\
f_3(n,2,3)<&1&1.1&1.2&1.3&1.4&1.6&1.8&2.1&2.4.\\
\end{array}$$
\vskip1mm

\ni From the bounds provided by this table and the properties of $f_3$ mentioned in the preceding paragraph we conclude that $f_3(n, d,3)<2.1$ for all $n\geqslant 7$, and $d\geqslant 2$. As ${D_k^{1/d}} <f_3(n,d,3)$, Proposition 2 implies that unless $k =\bQ$ (i.e., $d =1$),  $n\leqslant 6$.
\vskip1mm

We assert now that $n\leqslant 13$. To prove this, we can assume, in view of the result established in the preceding paragraph, that $k = \bQ$. By a direct computation we see that $f_1(14, 1, 1)<1$. Hence, $f_1(n,1,1)<1$ for all $n\geqslant 14$. As $D_{\bQ} =1$, from bound (15) we conclude that $n\leqslant 13$.  
  \vskip1mm
  
We will now assume that $d\geqslant 2$ and 
consider each of the possible cases  $2\leqslant n\leqslant 6$
separately.  
\vskip1mm

\ni $\bullet$ $n=6$: For $d\geqslant 2,$  $D_k^{1/d}<f_3(6,d,1)\leqslant f_3(6,2,1)<2.4.$
Therefore, by Proposition 2, $d=2$ and $D_k<6,$ which implies that  $k=\bQ({\sqrt 5})$ is the only possibility.
\vskip1mm

\ni $\bullet$ $n=5$:  For $d\geqslant 2,$ $D_k^{1/d}<f_3(5,d,1)\leqslant f_3(5,2,1)<2.9.$
Therefore, we infer from Proposition 2 that $d=2$ and $D_k<9$. So there are two possible real quadratic fields $k$, their discriminants are  $5$ and  $8.$  Both the fields have class number $1$, and  we use the
bound (15) to obtain $D_k^{1/2}<f_1(5,2,1)<2.8.$
So only $D_k=5$ can occur.
\vskip1mm

\ni $\bullet$ $n=4$: For $d\geqslant 3,$ $D_k^{1/d}<f_3(4,d,1)\leqslant f_3(4,3,1)<3.62$, and from Proposition 2 
we conclude that if $n=4$, then $d<3$. Let us assume that $d =2$. Then since $D_k^{1/2}< f_3(4,2, 1.1)<3.76$,  
$D_k<15$ and so the possible values of $D_k$ are $5, 8, 12$ or $13.$ The quadratic fields  with these $D_k$ have class number $1$. Now from 
bound (15) we obtain $D_k^{1/2}<f_1(4,2,1)<3.4.$
Hence, $D_k<12$, and only $D_k=5, 8$ can occur.

\vskip1mm

\ni $\bullet$ $n =3$: For $d\geqslant 4$, as $D_k^{1/d}<f_3(3,d,1)\leqslant f_3(3,4,1)<5.1$, from Proposition 2 we infer that if $n =3$, then $d<4$.  If $d = 3 =n$, $D_k<133$ from which we find that $D_k=49$ or $81.$
Now we consider the case where $d =2$ (and $n=3$). Since $D_k^{1/2}<f_3(3,2,1)<5.6$, $D_k<32$, and in this case the possible values of $D_k$ are $5, 8, 12, 13, 17,  21, 24, 28$ or $29.$ The quadratic fields with these discriminants have class number $1,$ and we use 
bound (15) to obtain $D_k^{1/2}<f_1(3,2,1)<4.52.$  Hence, 
$D_k<21$ and only $D_k=5, 8, 12, 13, 17$ can occur.
\vskip1mm

\ni $\bullet$ $n =2$: As $D_k^{1/d}<f_3(2,7,1)<9,$  
Proposition 2 implies that $d\leqslant 6$. 
\vskip1mm

$\vardiamondsuit$\:$n =2$ and $d=6$: As $D_k^{1/6}<f_3(2,6,1)<9$,  $D_k<531441.$  One can check from the table in  [1] that $h_k=1$ for all the five number fields satisfying this bound.  We now use bound (15) to obtain $D_k^{1/6}<f_1(2,6,1)<7.2$. But according to Proposition 2 there is no totally real number field $k$ for which this bound holds. 
\vskip1mm

 $\vardiamondsuit$\:$n=2$ and $d=5$: As $D_k^{1/5}<f_3(2,5,1)<9.3$,   $D_k<69569.$ 
Again, one can check from the table in [1] that there are five such number fields
and the class number of each of them is $1$.  Now we use 
bound (15) to obtain $D_k^{1/5}<f_1(2,5,1)<7.1.$  Hence, 
$D_k<18043.$  From [1] we find that $D_k=14641$ is the only possibility.
\vskip1mm

 $\vardiamondsuit$\:$n=2$ and $d=4$: As $D_k^{1/4}<f_3(2,4,0.9)<9.74$, 
$D_k<9000.$  According to [1], there are $45$ totally
real quartic number fields with  discriminant $<9000,$ all of them have 
class number $1.$  
We use 
bound (15) to obtain $D_k^{1/4}<f_1(2,4,1)<7.04.$  Hence, 
$D_k<2457.$  We find from [1] that
there are eight totally real quartic number fields $k$ with $D_k< 2457$. Their discriminants are  
$$725, 1125, 1600, 1957, 2000, 2048, 2225, 2304.$$
\vskip1mm

 $\vardiamondsuit$\:$n=2$ and $d=3$: As $D_k^{1/3}<f_3(2,3,0.8)<10.5$, 
$D_k<1158.$  From table B.4 of  [C] we find that there are altogether 
$31$ totally
real cubics satisfying this discriminant bound. Each of these fields have class number 
$1.$  We use 
bound (15) to obtain $D_k^{1/3}<f_1(2,3,1)<7$, which implies that  
$D_k<343.$  There are eight real cubic number fields
satisfying this bound. The values of $D_k$ are   
$$49, 81, 148, 169, 229, 257, 316, 321.$$ 
\vskip1mm

$\vardiamondsuit$\:$n=2$ and $d=2$: As $D_k^{1/2}<f_3(2,2,0.5)<12$, $D_k<144.$  From  table B.2 of totally real quadratic number
fields given in [C], we check that the class number of all
these fields are bounded from above by $2.$  Hence, 
$D_k^{1/2}<f_1(2,2,2)<7.3.$  So $D_k\leqslant 53.$  Among the real quadratic fields with $D_k\leqslant 53$, there is only one field whose class number is $2$, it is the field with $D_k=40$. All the rest have class number 1,  and from bound (15) we conclude that 
$D_k^{1/2}<f_1(2,2,1)<6.8$, i.e., $D_k<47$.  Therefore, the following is the list of the possible values of $D_k$:
$$5, 8, 12, 13, 17, 21, 24, 28, 29, 33, 37, 40, 41, 44.$$
\vskip3mm

To summarize, for $G$ of type $B_n$ or $C_n$, the  possible $n$, $d$ and $D_k$ are given in the following table.

$$\begin{array}{ccl}
n&\ \ \ \ d&\ \ \ \ D_k\\
2,\dots,13&\ \ \ \ 1&\ \ \ \ 1\\
6&\ \ \ \ 2&\ \ \ \ 5\\
5&\ \ \ \ 2&\ \ \ \ 5\\
4&\ \ \ \ 2&\ \ \ \ 5,8\\
3&\ \ \ \ 3&\ \ \ \ 49,81\\
3&\ \ \ \ 2&\ \ \ \ 5,8,12,13,17\\
2&\ \ \ \ 5&\ \ \ \ 14641\\
2&\ \ \ \ 4&\ \ \ \ 725, 1125, 1600, 1957, 2000, 2048, 2225, 2304\\
2&\ \ \ \ 3&\ \ \ \ 49, 81, 148, 169, 229, 257, 316, 321\\
2&\ \ \ \ 2&\ \ \ \ 5, 8, 12, 13, 17, 21, 24, 28, 29, 33, 37, 40, 41, 44.
\end{array}$$

\bs
\ni{\bf 4.3.} We will show that none of the possibilities listed in the above table actually give rise to an arithmetic fake compact hermitian symmetric space of type $B_n$ or $C_n$.  For this we recall first of all that $\overline{G}$, and so also $G$, is anisotropic over $k$ (1.5). Now we observe that if $G$ is a group of type $B_n$ ($n\geqslant 2$), then it is $k$-isotropic if and only if it is isotropic at all the real places of $k$ (this follows from the classical Hasse principle for quadratic forms which says that a quadratic form over $k$ is isotropic if and only if it is isotropic at every place of $k$, and the well-known fact that a quadratic form of dimension $>4$ is isotropic at every nonarchimedean place). Also, a $k$-group of type $C_n$ ($n\geqslant 2$) is $k$-isotropic if it is isotropic at all  the real places of $k$ (this is known, and follows, for example, from Proposition 7.1 of [PR]). These results imply that if $d =1$, i.e., if $k =\bQ$, then $G$ is isotropic, and so $d=1$ is not possible. 

\vskip1mm

Now let us take up the case where $d=2$, i.e., $k$ is a real quadratic field, and $n = 2, \,5$ or $6$. Then for any real place $v$ of $k$ where $G$ is isotropic, the complex dimension of the symmetric space of $G(k_v)$ is odd (recall from 1.4 that the complex dimension of the symmetric space of $G(k_v)$ is $2n-1$ if $G$ is of type $B_n$, and it is $n(n+1)/2$ if $G$ is of type $C_n$). But as the complex dimension of the hermitian symmetric space $X$ is even (since the orbifold Euler-Poincar\'e characteristic of $\Gamma$ is positive, see 1.3),  we conclude that $G$ must be isotropic at both the real places of $k$ (note that $G$ is anisotropic at a place $v$ of $k$ if and only if $G(k_v)$ is compact). From this observation we conclude that $G$ is $k$-isotropic also in case $d=2$, and $n =2,\, 5$ or $6$. Therefore these cases do not occur.

\vskip2mm

\ni{\bf 4.4.} To rule out the remaining cases listed in the table in 4.2, we compute  the value of $\cR$ in each case ($\cR$ as in (5)). The following table provides  the minimal monic polynomial defining $k$ and the values of $\zeta_k$ needed for the  computation of  $\cR$.  It turns out that in none of the remaining cases  the numerator of $\cR$ is a power of \,$2$ and Proposition 1 then eliminates these cases.

$$\begin{array}{ccccccccc}
n&d&k&D_k&\zeta_k(-1)&\zeta_k(-3)&\zeta_k(-5)&\zeta_k(-7)\\
4&2&x^2-5&5&1/30&1/60&67/630&361/120\\
4&2&x^2-2&8&1/12&11/120&361/252&24611/240.\\
\end{array}$$
\vskip2mm

$$\begin{array}{ccccccccc}
n&d&k&D_k&\zeta_k(-1)&\zeta_k(-3)&\zeta_k(-5)\\
3&3&x^3-x^2-2x+1&49&-1/21&79/210&-7393/63\\
3&3&x^3 - 3x - 1&81&-1/9&199/90&-50353/27\\
3&2&x^2 - 17&17&1/3&41/30&5791/63\\
3&2&x^2 - 13&13&1/6&29/60&33463/1638\\
3&2&x^2 - 3&12&1/6&23/60&1681/126\\
3&2&x^2 - 2&8&1/12&11/120&361/252\\
3&2& x^2 - 5&5&1/30&1/60&67/630.\\
\end{array}$$
\vskip2mm

$$\begin{array}{ccccccccccccc}
n&d&k&D_k&&\zeta_k(-1)&&\zeta_k(-3)\\
2&5&x^5-x^4-4x^3+3x^2+3x-1&14641&&-20/33&&1695622/165 \\
2&4&x^4 - 4x^2 + 1&2304&&1&&22011/10\\
2&4&x^4 - x^3 - 5x^2 + 2x + 4&2225&&4/5&&9202/5\\
2&4&x^4 - 4x^2 + 2&2048&&5/6&&87439/60\\
2&4&x^4 - 5x^2 + 5&2000&&2/3&&3793/3\\
2&4&x^4-4x^2-x+1&1957&&2/3&&3541/3\\
2&4&x^4 - 6x^2 + 4&1600&&7/15&&17347/30\\
2&4&x^4 - x^3 - 4x^2 + 4x + 1&1125&&4/15&&2522/15\\
2&4&x^4 - x^3 - 3x^2 + x + 1&725&&2/15&&541/15\\
2&3&x^3-x^2-4x+1&321&&-1&&555/2\\
2&3&x^3-x^2-4x+2&316&&-4/3&&874/3\\
2&3&x^3-x^2-4x+3&257&&-2/3&&1891/15\\
2&3&x^3 - 4x - 1&229&&-2/3&&1333/15\\
2&3&x^3 - x^2 - 4x -1&169&&-1/3&&11227/390\\
2&3&x^3 - x^2 - 3x + 1&148&&-1/3&&577/30\\
2&3&x^3 - 3x - 1&81&&-1/9&&199/90\\
2&3&x^3 - x^2 - 2x + 1&49&&1/21&&79/210.
\end{array}$$
\bs

 \centerline{\bf 5. $G$ of type $D_n$}
 \bs
 
We will consider hermitian symmetric spaces associated to Lie groups of type $D_n$, with $n\geqslant 4.$ The noncompact irreducible hermitian
symmetric spaces of these types are $\mathrm{SO}^*(2n)/\mathrm{U}(n)$ and $\mathrm{SO}(2, 2n-2)/\mathrm{S}(\mathrm{O}(2)\times \mathrm{O}(2n-2)).$  In the terminology of
\'Elie Cartan, these are hermitian symmetric spaces of types DIII and BDI respectively.

We note that any absolutely simple algebraic group $G$ over $\bQ$ of type $^1D_n$ or $^2D_n$, with $n\geqslant 4$, or a triality form of type $D_4$, whose $\R$-rank is at least $2$, is $\bQ$-isotropic (note that if $G$ is a triality form, then as at the unique real place of $\bQ$ the relative rank of $G$ is $2$, we see that in the Tits index of $G$ over $k_v$ the central vertex is distinguished for every place $v$ of $\bQ$, and then it follows from Proposition 7.1 of [PR] that $G$ is isotropic over $\bQ$), and hence, by Godement compactness criterion,  its arithmetic subgroups are non-cocompact in $G(\bR)$. Since we are only interested in compact hermitian locally symmetric spaces (and $\mathrm{SO}(2,2n-2)$ is of $\R$-rank 2, and for $n\geqslant 4$, $\R$-rank of $\mathrm{SO}^*(2n)$ is at least 2) {\it in this section the number field $k$ will be a nontrivial extension of} $\bQ$. 
\vskip2mm

\ni{\bf 5.1.}  The exponents of the Weyl group of $G$ of type $D_n$ are 
$1, 3, 5, \ldots, 2n-5, 2n-3,$ together with $n-1$ which has multiplicity two
if $n$ is even and multiplicity $1$ if $n$ is odd. The center of $G$
is of order $4$ and $\dim\,G = n(2n-1).$  Let $\cT$ be as in 2.10.
\vskip2mm

The following bounds for $[\Gamma:\Lambda]$ can be obtained from Propositions 0.12, 2.9, 5.1 and the considerations in 5.5 of [BP]. 
\vskip1mm

\ni{\bf Case (a):}  $n$ is even, and $G$ is of type $^1D_n$, i.e., it is of inner type. Then
\begin{equation}
[\Gamma:\Lambda] \leqslant h_{k,2}^22^{2(d+r-1+\#\cT)}.\end{equation}

\ni{\bf Case (b):} $n$ is even and $G$ is of type $^2D_n$. Then $\ell$ is a totally real quadratic extension of $k$ (see 1.5), and 
\begin{equation}
[\Gamma:\Lambda] \leqslant h_{\ell,2}2^{2(d+r+\#\cT)-1}D_{\ell}/D_k^2.\end{equation}

\ni{\bf Case (c):} $n$ is odd. Then  $G$ is of type $^2D_n$, $\ell$ is a totally complex quadratic
extension of $k$ (see 1.5), and  \begin{equation}
[\Gamma:\Lambda] \leqslant h_{\ell,4}2^{2(d+r+\#\cT)}.\end{equation}

\ni{\bf Case (d):} $n=4$, $G$ is a triality form of type $D_4$, $\ell$ is a totally real cubic extension of $k$ such that over the normal closure of $\ell/k$, $G$ is an inner form of a split group.
\begin{equation}
[\Gamma : \Lambda] \leqslant h_{\ell,2}2^{2(d+r+\#\cT)}D_{\ell}/D_k^3.\end{equation}
\vskip2mm

\begin{center}
{\bf Case (a)}
\end{center}

\ni{\bf 5.2.}  In this case, $n\,(\geqslant 4)$ is even, $G$ is of inner type, and (5) provides the following value of $\cR$:
\begin{equation*}
\cR :=2^{-dn}|\zeta_k(1-n)\prod_{j=1}^{n-1}\zeta_k(1-2j)|
\end{equation*}

Letting $$A(n)=\frac{(2\pi)^n}{(n-1)!}\cdot \prod_{j=1}^{n-1}\frac{(2\pi)^{2j}}{(2j-1)!}$$ and using the bounds (13) and (19) 
we obtain the following:

 \begin{equation}
D_k^{1/d} <
a_1(n,d,h_{k,2}):=\Big[ \{4A(n)\}^d\cdot \frac{h_{k,2}^2}4
\Big]^{\frac{2}{dn(2n-1)}}.
\end{equation}

Using the Brauer-Siegel bound (16) and lower bound for the regulator for
totally real field $k$ of degree $d$ recalled  in { 4.1} we obtain the following bound.
\begin{equation}
D_k^{1/d} <a_2(n,d,\delta)
:=\Big[\{\frac{[\Gamma((1+\delta)/2)\zeta(1+\delta)]^2}{(\pi)^{1+\delta}e^{0.92}}
A(n)
\}\cdot
\{25\delta(1+\delta)\}^{\frac{2}{d}} \Big]^{\frac{2}{(2n^2-n-2-2\delta)}} .
\end{equation}

The argument for the proof of  the following Lemma, which will be used in later sections as well, is the same as in  
the first paragraph of 4.2.
\begin{lemm}
Let $\delta\in [0.04,9].$  For fixed values of $n$ and $\delta,$ 
$a_2(n,d,\delta)$ decreases as $d$ increases.  Furthermore, for fixed values of $d$ and $\delta,$ if 
$n\geqslant 8,$  then $a_2(n+1,d,\delta)<\max(1, a_2(n,d,\delta)).$ 
\end{lemm}

We obtain by a direct computation
the following upper bound for
the value of $a_2(n,2,4)$ for small $n$.

$$\begin{array}{cccc}
n:&4&6&8 \\
a_2(n,2,4)<&10.7&3.33&2.13
\end{array}$$
From Proposition 2 we now infer that $k =\bQ$ for all even $n\geqslant 8$. But as $k\ne \bQ$, $n=4$ or $ 6$.
\vskip2mm

Consider first the case $n=6.$  For
$d\geqslant 2,$ $D_k^{1/d}< a_2(6,2,1)<3.2.$  Now using Proposition 2 we conclude that  $d=2$ and $D_k<11,$ hence, $D_k=5$ or $8.$  Since the class number of the corresponding fields is $1$, $D_k^{1/2}< a_1(6,2,1)<2.82.$  As $2.82^2<8,$
we conclude that $D_k=5.$ Then $k = \bQ(\sqrt {5})$ and for this field $\zeta_k(-1) = 1/30$, $\zeta_k(-3) = 1/60$, $\zeta_k(-5) = 67/630$, $\zeta_k(-7) = 361/120$ and $\zeta_k(-9)= 412751/1650$.  Using these values, we compute $\cR$ and find that its numerator is not a power of $2$, now Proposition 1 rules out the case $n =6$. 

\vskip2mm
Consider now $n=4.$  For $d\geqslant 4$,
  $D_k^{1/d}< a_2(4,4,1)<5.7.$   Therefore, $d\leqslant 4,$ and for $d=4,$ $D_k<1056.$  From the list of number
fields in [1] we find that the only possible value is $D_k=725$ and the class number of the corresponding number field 
is $1.$    Hence 
$D_k^{1/4}<  a_1(4,4,1)<4.9.$  According to Proposition 2 no such number field exists.

For $d=3,$ $D_k^{1/3}< a_2(4,3,1)<6.$  From the table of
totally real cubics in [1] 
we find that the class number of each of the four number fields satisfying the above bound is $1$.
Hence, 
$D_k^{1/3}< a_1(4,3,1)<5.$  So 
$D_k$ can only take one of the following two values, 
$$49,\:81.$$

For $d=2,$ $D_k^{1/2}< a_2(4,2,1)<6.7$; hence, $D_k< 45$.   From the list of
real quadratics in [C] we find for real quadratic $k$ with $D_k < 45$,  $h_k\leqslant 2.$  But then $D_k^{1/2}<  a_1(4,2,2)<5.$
We conclude that $D_k<25$ and then $h_k=1.$  It follows that $D_k^{1/2}< a_1(4,2,1)<4.73.$
We conclude that 
$D_k$ can only take one of the following values,
$$5,8,12,13,17,21.$$
Following is thus the list of possible totally real number fields 
$k$. 
$$\begin{array}{ccl}
n&d&D_k\\
4&3&49,81\\
4&2&5,8,12,13,17,21.
\end{array}
$$
In the following table, for each of these fields, we give the values of $\zeta_k$ required for the computation of $\cR$ for $n =4$

$$\begin{array}{ccccccccccc}
n&d&D_k&&\zeta_k(-1)&&\zeta_k(-3)&&\zeta_k(-5)\\
4&3&49&&-1/21&&79/210&&-7393/63\\
4&3&81&&-1/9&&199/90&&-50353/27\\
4&2&5&&1/30&&1/60&& 67/630 \\
4&2&8&&1/12&&11/120&&361/252\\ 
4&2&12&&1/6&&23/60&&1681/126\\
4&2&13&&1/6&&29/60&&33463/1638\\
4&2&17&&1/3&&41/30&&5791/63\\
4&2&21&&1/3&&77/30&&17971/63.
\end{array}
$$

\vskip2mm

Now computing $\cR$ we find that its numerator is not a power of $2$ and hence according to Proposition 1 none of the $k$ as above can give rise to  an arithmetic fake compact hermitian
space of type $D_n$ with $G$ of inner type.

\begin{center}
{\bf  Case (b)}
\end{center}

\ni{\bf 5.3.}  In this case, $n\,(\geqslant 4)$ is an even integer, $G$ is of type $^2D_n$, $\mathfrak{s} = 2n-1$, $s =4$, and $\ell$ is a totally real quadratic extension of $k$. The following value of $\cR$ is provided by (7):
\begin{equation*}
\cR =2^{-dn}|\zeta_{\ell|k}(1-n)\prod_{j=1}^{n-1}\zeta_k(1-2j)|.
\end{equation*}

Letting $$A(n)=\frac{(2\pi)^n}{(n-1)!}\cdot \prod_{j=1}^{n-1}\frac{(2\pi)^{2j}}{(2j-1)!}$$ and using the bounds (13) and (20), we obtain the following bounds:
 \begin{equation}
D_k^{1/d} <
b_1(n,d,h_{\ell,2}):=\Big[ \{4A(n)\}^d\cdot \frac{h_{\ell,2}}2
\Big]^{\frac{2}{dn(2n-1)}},
\end{equation}
\begin{equation}
D_k^{1/d} < b_2(n,d,\delta):=\Big[\{\frac{[\Gamma((1+\delta)/2)\zeta(1+\delta)]^2}{(\pi)^{1+\delta}e^{0.92}}
A(n)\}\cdot
\{25\delta(1+\delta)\}^{\frac{1}{d}} \Big]^{\frac{2}{(2n^2-n-2-2\delta)}},
\end{equation}
\begin{equation}
D_\ell^{1/2d}<{\mathfrak t}_1(n,d,D_k,h_{\ell,2})\\
:=\Big[ 2^{2d-1}A(n)^d h_{\ell,2}D_k^{\frac{n(5-2n)-6}{2}}\Big]^{\frac{1}{d(2n-3)}},
\end{equation}
\begin{eqnarray}
D_\ell^{1/2d}&<&{\mathfrak t}_2(n,d,D_k,R_\ell/w_\ell,\delta)\\
\nonumber&:=&\Big[ \frac{\delta(1+\delta)}{2R_\ell/w_\ell}D_k^{\frac{n(5-2n)-6}{2}}\{A(n)\frac{[\Gamma({(1+\delta)}/2)\zeta(1+\delta)]^2}{\pi^{1+\delta}}\}^d\Big]^{\frac{1}{d(2n-4-\delta)}},
\end{eqnarray}
\begin{equation}
D_\ell^{1/2d}<{\mathfrak t}_3(n,d,D_k,\delta):=\Big[25\delta(1+\delta)D_k^{\frac{n(5-2n)-6}{2}}\{A(n)\frac{[\Gamma({(1+\delta)}/2)\zeta(1+\delta)]^2}{\pi^{1+\delta}e^{0.92}}\}^d\Big]^{\frac{1}{d(2n-4-\delta)}}.
\end{equation}
To obtain the above bounds we have used $D_{\ell}\geqslant D_k^2$, the Brauer-Siegel bound (16) for totally real fields, and the bound for the regulator due to Zimmert given in 
{4.1}.  

\bs
\ni{\bf 5.4.}  We obtain by a direct computation
the following upper bound for
the value of $b_2(n,2,4)$ for small $n$.

$$\begin{array}{cccc}
n:&4&6&8\\
b_2(n,2,4)<&7.6&3&2.1.
\end{array}$$

   From Proposition 2, and Lemma 1, where in the latter the function $a_2(n,d,\delta)$ is replaced by $b_2(n,d,\delta),$ 
we conclude that $k=\bQ$ for all even $n\geqslant8.$ But $k\ne \bQ$, and hence $n\leqslant 6$.
\vskip2mm

Consider now $n=6.$  For $d\geqslant3,$ $D_k^{1/d}< b_2(6,2,2)<3.$  Therefore, $d=2$ and $D_k=5, 8$ are the only possibilities. Let us take up first the case where $D_k=8.$  As   
$D_\ell^{1/4}\leqslant {\mathfrak t}_3(6,2,8,2)<3.2$,  Proposition 2 rules out this case. 
Consider now the case where $D_k=5$, i.e., $k =\bQ(\sqrt{5})$. The following argument involving Hilbert class fields
will be used repeatedly.  As 
$D_\ell^{1/4}\leqslant {\mathfrak t}_3(6,2,5,1)<7.3.$ The Hilbert class field of $\ell$  is a totally real number field (since $\ell$ is totally real) of degree $h_{\ell}$ over $\ell$ (hence of degree 
$4h_{\ell}$ over $\bQ$), and its root discriminant equals $D_{\ell}^{1/4}$ which is $<7.3.$   On the other hand,
 according to Proposition 2 $M_r(6)> 8.18.$  So we conclude that $4h_\ell<6$ and,
$h_\ell\leqslant \ll5/4\lr=1.$  Where here, and in the sequel,  we use $\ll x \lr$ to denote the integral part of $x.$ 
It follows that
$D_\ell^{1/4}\leqslant {\mathfrak t}_1(6,2,5,1)<5.5.$  From [1] we see that there is  only one number field  $\ell$ containing $k =\bQ(\sqrt{5})$ with this root discriminant bound. For this $\ell$, 
$D_\ell=725$, and $\zeta_{\ell |k}(-5) = 2164$.  Now using this value of $\zeta_{\ell |k}(-5)$ and the values of $\zeta_k$ (for $k =\bQ(\sqrt{5})$) given in  5.2 we compute the value of $\cR$ and find that its numerator is not a power of $2$. So Proposition 1 rules out $n =6$ with $D_k =5$.

\vskip2mm

Consider now $n=4.$   As $D_k^{1/d}<b_2(4,4,2)<5.17,$ 
Proposition 2 implies that $d\leqslant 3$.  

For $d=3,$  we know from Proposition 2 that $D_k\geqslant 49.$ Hence,  
$D_\ell^{1/6}\leqslant {\mathfrak t}_3(4,3,49,1)<17.$  From Table IV of [Mart], $M_r(14)>17.$ 
So by considering the Hilbert class field of $\ell$, we obtain  
$h_\ell\leqslant \ll13/6\lr =2.$ 
It follows that
$D_\ell^{1/6}\leqslant {\mathfrak t}_1(4,3,49,2)<8.6.$  But according to Proposition 2, $M_r(7)>11.05.$ Hence, 
$h_\ell\leqslant \ll8/6\lr=1.$  This in turn implies that $D_\ell^{1/6}< {\mathfrak t}_1(4,3,49,1)<8.2$, and therefore, $D_{\ell}< 304007$.  
It is seen from table t66.001 in [1] that there is only one totally real number field  $\ell$ of degree $6$ for which this bound holds. For this $\ell$, $D_\ell=300125.$
Hence the only possibility for $d=3$ is 
$(D_k, D_\ell)=(49, 300125).$
\vskip1mm

Let us assume now that $d=2$. For a real quadratic field $k$, either $D_k = 5\: {\mathrm {or}}\:8$ or $D_k\geqslant 12$. 
Consider first
the quadratic fields $k$ with $D_k\geqslant 8.$ Since 
$D_\ell^{1/4}< {\mathfrak t}_3(4,2,8,0.5)<39.2.$ From Table IV of [Mart] we find that $M_r(80)>39.4.$  Hence, 
by considering the Hilbert class field of $\ell$, we infer that  
$h_\ell\leqslant \ll79/4\lr=19.$  Hence $h_{\ell,2}\leqslant 16.$
It follows that
$D_\ell^{1/4}\leqslant {\mathfrak t}_1(4,2,8,16)<16.79.$ As $M_r(14)>17$, by considering the Hilbert class field of $\ell$, we conclude that  
$h_\ell\leqslant \ll13/4\lr=3.$  So $h_{\ell,2}\leqslant 2.$ But then 
$D_\ell^{1/4}\leqslant {\mathfrak t}_1(4,2,8,2)<13.637$ and hence, $D_\ell\leqslant 34584.$
\vskip1mm

Let us now consider real quadratic fields $k$ with $D_k\geqslant 12$.  The discussion in the preceding paragraph 
implies that $h_{\ell,2}\leqslant 2.$  As 
$D_\ell^{1/4}\leqslant {\mathfrak t}_1(4,2,12,2)<9.47.$ Proposition 2 gives that  $M_r(8)>11.38.$  Hence,
by considering the Hilbert class field of $\ell$, we conclude that  
$h_\ell\leqslant \ll7/4\lr=1.$  But then as 
$D_\ell^{1/4}\leqslant {\mathfrak t}_1(4,2,12,1)<8.834$, so $D_\ell\leqslant 6090.$

 From t44.001 again, we check that
there are only $24$ such totally real quartics, with
$D_\ell$ given below:
$$725, 1125, 1600, 1957, 2000, 2048, 2225, 2304, 2525, 2624, 2777, 3600,$$ $$3981, 4205, 4225, 4352,4400, 4525, 4752, 4913, 5125, 5225, 5725, 5744.$$
Furthermore, with our assumption that $D_k\geqslant 12,$ we know that $h_{\ell} =1$, and as  
$D_k^{1/2}\leqslant b_1(4,2,1)<4.84$, $D_k$ can only be one of  $12, 13, 17, 21$.  
Since $D_\ell$ is an  integral multiple of $D_k^2,$ we check 
easily that for $D_k\geqslant 12,$ the only possible values for  $(D_k, D_\ell)$ are 
$(17,4913), (13, 4225), (12, 2304),$ $(12, 3600)$ and $(12,4752).$
\vskip1mm

Consider now the case  $D_k=5$, i.e., $k =\bQ(\sqrt{5})$, and $\ell$ is a totally real number field of degree $4$ containing $\bQ(\sqrt{5})$.  We will show that $D_\ell^{1/4}\leqslant 55.$ 
Assume to the contrary that $D_\ell^{1/4}>55.$   We will first prove that $R_{\ell}\geqslant 1.64$.  For this we shall use some results of [F], \S3. In the following paragraph all unexplained notation are from [F], \S 3, in which $k$ has been replaced by $\ell$.  

Recall that the image of the group of units of $\ell$ under the logarithmic embedding $\ell-\{0\}\rightarrow \bR^4$ 
forms a lattice $\Lambda_\ell$ of rank $3.$  Let $0<m_\ell(\varepsilon_1)\leqslant m_\ell(\varepsilon_2)\leqslant m_\ell(\varepsilon_3)$
be the successive minima of the Euclidean abasolute value on $\Lambda_{\ell}.$  Consider first the case where $\bQ(\varepsilon_1)=\ell.$  In this case, using Remak's estimate
as stated in (3.15) of [F], we see that 
the following lower bound for the regulator of $\ell$ holds:
$$R_\ell\geqslant \Big(\frac{\log D_\ell-4\log 4}{40^{1/2}}\Big)^3>4.5.$$

Let us assume now that $\bQ(\varepsilon_1)$ is a proper subfield of $\ell.$  Then $\bQ(\varepsilon_1)$
is a real quadratic field.  Among such fields, $\bQ(\sqrt{5})$ has the smallest regulator (see the corollary in \S 3 of [Z]).  Hence, 
the smallest fundamental unit $\varepsilon_1$ can be taken to be  $\frac{1+\sqrt{5}}{2}$. Then 
$m_\ell(\varepsilon_1)=2\log(\frac{1+\sqrt{5}}2)$.  So, 
$m_\ell(\varepsilon_2)\geqslant 2\log(\frac{1+\sqrt{5}}2).$ From a result of Remak and Friedman, cf.\:(3.2) of [F], we know that
$$m_\ell(\varepsilon_3)\geqslant 2(\frac14\log|D_\ell|-\frac12\log 5-\log2)>2(\log(55)-\log(2\sqrt5)).$$
(Note that $A(\ell/k)=(\frac23(8-2))^{1/2}={2}$ in the notation of [F], page 611.)
Hence from the bound (3.12) of [F] we obtain the following:
\begin{eqnarray*}
R_\ell&\geqslant&\frac1{2\sqrt2}\prod_{i=1}^3m_\ell(\varepsilon_i)\geqslant \sqrt2\Big(\log(\frac{1+\sqrt{5}}2)\Big)^2m_\ell(\varepsilon_3)\\
&>&2\sqrt2\Big(\log(\frac{1+\sqrt{5}}2)\Big)^2\Big(\log(55)-\log(2\sqrt5)\Big)
>1.64.\\
\end{eqnarray*}
This proves our assertion about $R_{\ell}$. Now since $w_\ell=2,$ we conclude that 
$$D_\ell^{1/4}< {\mathfrak t}_2(4,2,5,1.64/2,0.5)<55,$$ contradicting the assumption
that $D_\ell^{1/4}>55.$  Thus  we have proved that $D_\ell^{1/4}\leqslant 55.$
\vskip1mm

We find, using Table IV of [Mart], that $M_r(800)>55.$  As $D_{\ell}^{1/4}\leqslant 55$, 
by considering Hilbert class field of $\ell$, we conclude that  
$h_\ell\leqslant \ll799/4\lr=199.$  Hence, $h_{\ell, 2}\leqslant 128$. Then 
$D_\ell^{1/4}< {\mathfrak t}_1(4,2,5,128)<31.6.$ From Table IV of [Mart] we see that $M_r(41)>31.7.$  Hence,
again by considering the Hilbert class field of $\ell$ we conclude that  
$h_\ell\leqslant \ll40/4\lr=10$, so  $h_{\ell,2}\leqslant 8.$
But then 
$D_\ell^{1/2d}\leqslant {\mathfrak t}_1(4,2,5,8)<24.$
Again, from Table IV of [Mart] we see that  $M_r(24)>24$.   
By considering the Hilbert class field of $\ell$, we infer that 
$h_\ell\leqslant \ll23/4\lr=5.$  Hence, $h_{\ell,2}\leqslant 4.$ But then 
$D_\ell^{1/4}< {\mathfrak t}_1(4,2,5,4)<22.32$, and so $D_\ell< 248186.$  
From the tables t44001-t44003 of [1] we find that for $D_\ell\leqslant 248186,$
$h_\ell\leqslant 3$, and so 
$h_{\ell,2}\leqslant 2.$
It then follows that $D_\ell\leqslant \ll{\mathfrak t}_1(4,2,5,2)^4\lr\leqslant 187789.$

Here is the list of all the possibilities:
$$
\begin{array}{ccll}
n&d&D_k&D_\ell\\
4&3&49&300125\\
4&2&17&4913\\
4&2&13&4225\\
4&2&12&2304, 3600, 4752\\
4&2&8&\leqslant 34584\\
4&2&5 &\leqslant 187789.
\end{array}
$$

\ni{\bf 5.5.}  Malle has provided us the list of pairs $(k,\ell)$ satisfying the above 
constraints.   The values of $\zeta_k$ and $\zeta_{\ell |k}$ required to compute $\cR$ for each of the possible pairs $(k,\ell )$, with $D_k\geqslant 12$, have been  tabulated below.

$$\begin{array}{ccccccccccccc}
n&d&D_k&&D_\ell&&\zeta_k(-1)&&\zeta_k(-3)&&\zeta_k(-5)&&\zeta_{\ell|k}(-3)\\
4&3&49&&300125&&-1/21&&79/210&&-7393/63&& 8202104\\
4&2&17&&4913&&1/3&& 41/30&&5791/63&&366280/17\\
4&2&13&&4225&&1/6&&29/60&&33463/1638&&35936\\
4&2&12&&2304&&1/6&&23/60&&1681/126&& 5742\\
4&2&12&&3600&&1/6&&23/60&&1681/126&&25776\\
4&2&12&&4752&&1/6&&23/60&&1681/126&&68944.
\end{array}$$

Using the values of $\zeta_k$ and $\zeta_{\ell |k}$ given above, we can compute $\cR$. We see that its numerator is 
not a power of $2$ for any of the  above $(k,\ell )$, and Proposition 1 rules out all these pairs.
\vskip2mm

For $k=\bQ(\sqrt{2})$, for which $D_k=8$, there are $32$ number fields $\ell$ containing $k$ and with $D_{\ell}\leqslant 34584$. 
For $k=\bQ(\sqrt{5}),$ for which $D_k=5,$ there are $363$ number fields $\ell$ containing $k$ and with $D_{\ell}\leqslant 187789$.   In each of these $32+363$ cases, we have computed $\cR$  (interested readers my write to either of the authors to
obtain the values).  The numerator of $\cR$ in none of the cases is a power of $2$.  Proposition 1 thus eliminates Case (b).

\bs

\begin{center}
{\bf Case (c)}
\end{center}
\bs

\ni{\bf 5.6.}  In Case (c), $n$ is odd and $G$ is of type $^2D_n$, $\mathfrak{s} = 2n-1$, $s = 4$, $\ell$ is a totally complex quadratic extension of totally real $k\: (k \ne \bQ)$.  Equation (7) provides the following value of $\cR$:
\begin{equation*}
\cR =2^{-dn}|\zeta_{\ell|k}(1-n)\prod_{j=1}^{n-1}\zeta_k(1-2j)|.
\end{equation*}

Letting $$A(n)={\frac{(2\pi)^n}{(n-1)!}}\cdot \prod_{j=1}^{n-1}\frac{(2\pi)^{2j}}{(2j-1)!},$$ from bounds (13) and (21), 
using the following bound provided by the  Brauer-Siegel Theorem for a totally complex 
number field $\ell$ of degree $2d,$ 
\begin{equation}
 h_\ell R_\ell\leqslant w_\ell \delta(1+\delta)\Gamma(1+\delta)^d((2\pi)^{-2d}D_\ell)^{(1+\delta)/2}
\zeta_\ell(1+\delta),
\end{equation}
where  $\delta>0$, $h_\ell$ is the class number and $R_\ell$ is the regulator of $\ell$, 
and
$w_\ell$ is the order of the finite group of roots of unity contained 
in $\ell,$ and the bound 
$R_\ell\geqslant 0.02 w_\ell\,e^{0.1d}$ due to R.\,Zimmert [Z], we obtain the following bounds:

\begin{equation}
D_k^{1/d} <
c_1(n,d,h_{\ell,4}):=\Big[ \{4A(n)\}^d h_{\ell,4}
\Big]^{\frac{2}{dn(2n-1)}},\end{equation}
\begin{equation}
D_k^{1/d} <c_2(n,d,\delta):=\Big[\{4A(n)\frac{\Gamma(1+\delta)\zeta(1+\delta)^2}{(2\pi)^{1+\delta}e^{0.1}}\}\cdot
\{50\delta(1+\delta)\}^{\frac{1}{d}} \Big]^{\frac{2}{2n^2-n-2-2\delta}},\end{equation}
\begin{equation}
D_{\ell}/D_k^2<\mathfrak{t}(n,d,D_k,h_{\ell,4}):=\Big(4^dA(n)^dh_{\ell,4}\Big)^{\frac{2}{2n-1}}D_k^{-n},\end{equation}
\begin{equation}
D_\ell^{1/2d}<{\mathfrak u}_1(n,d,D_k,h_{\ell,4}):=\Big[\Big(4^dA(n)^dh_{\ell,4}\Big)^{\frac{2}{2n-1}}D_k^{2-n}\Big]^{\frac{1}{2d}},\end{equation}
\begin{eqnarray}D_\ell^{1/2d}&<&{\mathfrak u}_2(n,d,D_k,R_\ell/w_\ell,\delta)\\
\nonumber&:=&\Big[ \frac{\delta(1+\delta)}{R_\ell/w_\ell}D_k^{\frac{n(5-2n)-2}{2}}\{4A(n)\frac{\Gamma(1+\delta)\zeta(1+\delta)^2}{(2\pi)^{1+\delta}}\}^d\Big]^{\frac{1}{d(2n-2-\delta)}},\end{eqnarray}
\begin{equation}
D_\ell^{1/2d}<{\mathfrak u}_3(n,d,D_k,\delta):=\Big[50\delta(1+\delta)D_k^{\frac{n(5-2n)-2}{2}}\{4A(n)\frac{\Gamma(1+\delta)\zeta(1+\delta)^2}{(2\pi)^{1+\delta}e^{0.1}}\}^d\Big]^{\frac{1}{d(2n-2-\delta)}}.\end{equation}

\bs
\ni{\bf 5.7.} We obtain by a direct computation
the following upper bound for $c_2(n,2,3)$ for small $n$.

$$\begin{array}{cccc}
n:&5&7&9\\
c_2(n,2,2.6)<&4.2&2.5&1.78.\\
\end{array}$$
\vskip1mm

It is obvious that the conclusion of Lemma 1 holds with the function $a_2(n,d,\delta)$ replaced by  $c_2(n,d,\delta)$. Also, for fixed $d$ and $\delta$, $c_2((n,d,\delta)$ clearly decreases as $n$ increases. As $k\ne \bQ$,  using Proposition 2 we conclude that  $n\leqslant 7.$
\vskip1mm

Consider now $n=7.$  For $d\geqslant 2$, since   
$D_k^{1/d}\leqslant c_2(7,2,2)<2.5,$ Proposition 1 implies that $d=2$
and $D_k=5$ is the only possibilty. 
But if $D_k=5,$ $D_\ell^{1/4}\leqslant {\mathfrak u}_3(7,2,5,1.5)\leqslant 3.2,$ which according to Proposition 2 is not possible. 

\vskip2mm

\ni{\bf 5.8.}  Consider now the case 
 $n=5.$  As $c_2(5,3,1)<4,$ Proposition 2 implies that $d\leqslant 3$. If $d =2$, then $D_k < 4^2 = 16$, and hence,  $D_k=5,8,12$ or $13.$ On the other hand,  if $d = 3$, then $D_k <4^3 = 64$, and $D_k=49$ is the only possibility.  But then  
$D_\ell^{1/6}\leqslant {\mathfrak u}_3(5,3,49,1)<4.4.$ According to Proposition 2 there does not exist a totally complex  
$\ell$ of degree $6$ satisfying this bound for $D_{\ell}$.  We conclude therefore that $d=2.$
\vskip1mm

Let now $d =2$. If $D_k\geqslant 8$, as $D_{\ell}^{1/4}< \mathfrak{u}_1(5,2,8,2)< 5.52$, and according to Proposition 2, $M_c(8) > 5.78$, we conclude using the Hilbert class field of $\ell$ that  $h_\ell\leqslant \ll7/4\lr=1.$
\vskip1mm

Let us consider the case $D_k=13$.
As ${\mathfrak t}(5,2,13,1)<1.1,$ $D_\ell=169.$  However, there
is no totally complex quartic field with discriminant $169.$  Hence $D_k$ cannot
be $13.$

Let us now assume that $D_k=12$. 
As ${\mathfrak t}(5,2,12,1)<1.7,$ 
the only possibility for $\ell$ is $D_\ell=144$.

Suppose $D_k=8.$ Then $k =\bQ(\sqrt{2})$. As ${\mathfrak t}(5,2,8,1)<
12.36$, $D_\ell\leqslant 12\cdot 8^2\leqslant768.$   From the list of totally complex quartics  in table t40.001 of [1] we see that  those which contain $\bQ(\sqrt{2})$, and have discriminant in the above range, have discriminant in $ \{256, 320, 512, 576\}$. (Note that there are two totally complex quartics with discriminant $576$, but only one of them contains $\bQ(\sqrt{2})$. Only the one containing $\bQ(\sqrt{2})$ is of interest to us.)  

Now let us assume that $D_k =5$. Then $k = \bQ(\sqrt{5})$. As     
$D_\ell^{1/4}\leqslant {\mathfrak u}_3(5,2,5,0.7)<12.4$, $M_c(36)>12.5$ (Table IV of [Mart]),  
by considering the Hilbert class field of $\ell$, we infer that $h_\ell\leqslant \lfloor35/4\rfloor=8.$ But then as 
$D_\ell^{1/4}\leqslant {\mathfrak u}_1(5,2,5,8)<
8.5$, and $M_c(16)>8.7$ (Table IV of [Mart]), 
by again considering the Hilbert class field of $\ell$, we conclude that $h_\ell\leqslant \lfloor15/4\rfloor=3.$
It follows that $h_{\ell,4}\leqslant 2$ and hence, 
$D_\ell^{1/4}<  {\mathfrak u}_1(5,2,5,2)<7.85$. Therefore,  $D_\ell \leqslant 3797$. Moreover, $D_{\ell}$ is a multiple of  $D_k^2 =25$, so $D_{\ell}\leqslant 3775$.   According to the table in [1],  The discriminant of the totally complex quartic fields $\ell$ with $D_{\ell}\leqslant 3775$, and which contain $\bQ(\sqrt{5})$, is one of the following: 
$$125,225,400,1025,1225,1525,1600,2725,3025, 3625, 3725.$$
The class number of $\ell$ with $D_\ell=3725$ is $1$ according to [Mart].  But ${\mathfrak t}(5,2,5,1)<130$ and hence
$D_\ell\leqslant 5^2\cdot129< 3225.$  Hence $\ell$ with $D_{\ell} =3725$ can be excluded.

\vskip2mm

\ni{\bf 5.9.}   Among the pairs of $(D_k, D_\ell)$ obtained above, only some
of them can be discriminants of number fields $k$ and $\ell$ such that
$\ell$ is a totally complex quadratic extension of $k.$   We eliminate the rest.
In conclusion, here are all the possibilities in 
Case (c): $n =5$, $d=2$, $\ell$ is a totally complex quadratic extension of a 
real quadratic number field $k$, and    
$$
\begin{array}{cccl}
n&d&D_k&D_\ell\\
5&2&12&144\\
5&2&8&256,576\\
5&2&5&125,225,400,1025,1225,1525,1600,2725,3025, 3625.
\end{array}
$$

Using the values of $\zeta_k$ and $\zeta_{\ell | k}$ given in the following table for the pairs $(k,\ell)$ listed above, we computed $\cR$. Its numerator for none of the pairs $(k,\ell)$ turned out to be a power of $2$.  So by Proposition 1,  Case (c) does not  give rise to any arithmetic fake compact hermitian symmetric spaces.

$$\begin{array}{cccccccccccccc}
D_k&D_\ell&&\zeta_k(-1)&&\zeta_k(-3)&&\zeta_k(-5)&&\zeta_k(-7)&&\zeta_{\ell | k}(-4)\\
12&144&&1/6&&23/60&&1681/126&&257543/120&&5/3\\
8&256&&1/12&&11/120&&361/252&&24611/240&&285/2\\
8&576&&1/12&&11/120&&361/252&&24611/240&&15940/3.\\
\end{array}$$

For $k=\bQ(\sqrt{5}),$\,  $\zeta_k(-1)=1/30,\, \zeta_k(-3)=1/60,\, \zeta_k(-5)=67/630, \,\zeta_k(-7)=361/120$, and  

$$\begin{array}{cccccc}
D_\ell:&125&225&400&1025&1225\\
\zeta_{\ell | k}(-4):&1172/25&1984/3&8805&608320&1355904\\
\end{array}$$

$$\begin{array}{cccccccc}
D_\ell:&1525&1600&2725&3025&3625\\
\zeta_{\ell | k}(-4):&3628740&4505394&49421124
&872059200/11&178910784.
\end{array}$$
\bs

\begin{center}
{\bf Case (d)}
\end{center}

\bs
\ni{\bf 5.10.} We shall finally consider triality forms of type $D_4$. So assume now that $G$ is a triality form over a totally real number field $k\ne \bQ$. For such a $G$, $\mathfrak{s} =7$, $s =4$, $\dim\, G = 28$, and $\ell$ is a totally real cubic extension of $k$ such that over the normal closure of $\ell/k$, $G$ is an inner form of a split group. 

\vskip2mm

The exponents of the Weyl group of $G$ are $1$, $3$, $3$, and $5$ ($3$ has multiplicity $2$). 
The value of $\cR$ in this case, provided by equation (8), is 
$$\cR=2^{-4d}|\zeta_k(-1)\zeta_{\ell|k}(-3)\zeta_{k}(-5)|.$$

Letting $A=(2\pi)^{16}/4320$ and using bounds (13), (16) and (22), and Zimmert's lower bound for the regulator,  we conclude that
\begin{equation}
D_k^{1/d} <d_1(d,h_{\ell,2}):=\Big[ (4A)^d h_{\ell,2}
\Big]^{1/14d},\end{equation}
\begin{equation}
D_k^{1/d} <d_2(d,\delta):=[(50\delta(1+\delta))\Big(A\frac {\{\zeta(1+\delta)\Gamma(\frac{1+\delta}2)\}^3}{2e^{1.38}
\pi^{\frac{3}{2}(1+\delta)}}\Big)^d]^{\frac2{d(25-3\delta)}}, \end{equation}
\begin{equation}
D_\ell^{1/3d}<\mathfrak{z}_1(d,D_k,h_{\ell,2}):=((4A)^dD_k^{-\frac{13}2}h_{\ell,2})^{\frac2{15d}},\end{equation}
\begin{equation}
D_\ell^{1/3d}<\mathfrak{z}_2(d,D_k,\delta):=
[50\delta(1+\delta)D_k^{-\frac{13}2}\Big(A\frac{\{\zeta(1+\delta)\Gamma(\frac{1+\delta}2)\}^3}{2e^{1.38}\pi^{\frac32(1+\delta)}}\Big)^d]^{\frac2{3d(4-\delta)}}.\end{equation}

Note that for a fixed value of $\delta\geqslant 0.02$, all the expressions on the right hand side
of the above bounds are decreasing in $d$.  By a 
direct computation we find that 
$D_k^{1/d}< d_2(4,1.6)<5.03.$ Using Proposition 2 we conclude from this that $d<4$.
\vskip1mm

Consider now $d=3.$  As the smallest discriminant of a 
totally real cubic is $49,$ and  
$D_\ell^{1/9}<\mathfrak{z}_2(3,49,1)<10$. But 
$M_r(9)>11.8$ (see Table IV in
[Mart]).  Hence $d$ cannot be $3$.
\vskip1mm

Consider now $d=2.$  As the smallest discriminant of a 
totally real quadratic field is $5,$ and  
$D_\ell^{1/6}<\mathfrak{z}_2(2, 5, 0.7)<22.2.$   But 
$M_r(21)>22.3$ (Table IV in [Mart]).  So by considering the Hilbert class field of $\ell$, we conclude that 
$h_\ell\leqslant\ll20/6\lr=3$,  and hence, $h_{\ell,2}\leqslant 2.$ 
But  then 
$D_\ell^{1/6}<\mathfrak{z}_1(2,5,2)<10.4.$
According to Proposition 2, 
$M_r(8)>11.38$, so considering again the  Hilbert class field of $\ell$, we infer that 
$h_\ell\leqslant\ll7/6\lr=1.$  Now since 
$D_\ell^{1/6}<\mathfrak{z}_1(2,5,1)<9.896$, it follows that $D_\ell\leqslant \ll9.896^6\lr=939200.$
\vskip1mm

Suppose $D_k\geqslant 8,$ then  
$D_\ell^{1/6}< \mathfrak{z}_1(2,8,1)<8.1,$ which is smaller than
the lower bound for $M_r(6)$ given by Proposition 2.
Hence $D_k$ can only be  $5$.
\vskip1mm

For $D_k=5,$ i.e., $k =\bQ(\sqrt{5})$, we find from table t66.001 of [1]
that there are $11$ totally real sextics with discriminant
bounded as above.  For $\ell$ to be an extension of degree three
of $k,$ it is necessary that $D_\ell/D_k^3$ is an integer.  Going
through the list of the $11$ sextics, we are left with four possibilities for $D_\ell$, these are $300125, 485125, 722000$ and
$820125.$  Among these four sextics, 
only the  one with $D_{\ell} = 300125$ contains $\bQ(\sqrt {5})$ as a subfield.  This $\ell$ is given by $x^6-x^5-7x^4+2x^3+7x^2-2x-1$.  The values of  $\zeta_k(-1)$, $\zeta_{\ell |k}(-3)$, $\zeta_k(-5)$ and $\cR$ are given below.

$$\begin{array}{ccccccccc}
\zeta_k(-1)&&\zeta_k(-5)&\zeta_{\ell|k}(-3)&\cR\\
1/30&&67/630&1295932432/7&5426717059/33075.\\
\end{array}$$
As the numerator of $\cR$ is not a power of $2$, from Proposition 1 we conclude that 
arithmetic fake compact hermitian
symmetric spaces of type $D_4$ cannot arise from  triality forms.

\bs
\ni{\bf 5.11.}
In conclusion, there does not exist  an arithmetic fake compact hermitian symmetric space
of type $D_n$, \,$n\geqslant 4.$

\vskip6mm 

\begin{center} {\bf 6. $G$ of type $^2E_6$ }  
\end{center}

\bs
\ni{\bf 6.1.} In this section $G$ is of type $^2E_6$. Its dimension is $78$ and the complex dimension of 
the symmetric space of $\cG =\prod_{j =1}^r G(k_{v_j})$ is $16r$.  The exponents of the Weyl group of $G$ are $1,\,4,\,5,\,7,\,8$ and $11$, $\mathfrak{s} = 26$, and $s = 3$. Let $\cT$ be as in 2.10. The bound (13) in the present case is 
\begin{equation}
(D_kD_{\ell})^{13}\Big(\frac{4!5!7!8!11!}{(2\pi)^{42}}\Big)^d3^{\#\cT} < [\Gamma : \Lambda]/3^r.\end{equation}
\vskip1mm

\ni {\bf 6.2.} The center $C$ of $G$ is  $k$-isomorphic to the kernel of the norm map $N_{\ell/k} : R_{\ell /k}(\mu_3) \to \mu_3$. As this map is onto, the Galois cohomology group $H^1(k, C)$ is isomorphic to the kernel  of the homomorphism $\ell^{\times}/{\ell^{\times}}^3\to k^{\times}/{k^{\times}}^3$ induced by the norm map. We shall denote this kernel by $(\ell^{\times}/{\ell^{\times}}^3)_{\bullet}$. 

By Dirichlet's unit theorem, $U_k\cong \{\pm 1\}\times 
{\bZ}^{d-1}$, and $U_{\ell}\cong \mu(\ell)\times {\bZ}^{d-1}$, where
$\mu(\ell)$ is the finite cyclic group of roots of unity in $\ell$. Hence, 
$U_k/U_k^3\cong (\bZ/3\bZ)^{d-1}$, and $U_{\ell}/U_{\ell}^3\cong
\mu(\ell)_3\times (\bZ/3\bZ)^{d-1}$, where $\mu(\ell)_3$ is the group
of cube-roots of unity in $\ell$. Now we observe that
$N_{\ell/k}(U_{\ell})\supset N_{\ell/k}(U_k) = U_k^2$, which implies that 
the homomorphism  $U_{\ell}/U_{\ell}^3\rightarrow 
U_k/U_k^3$, induced by the norm map, is onto. Therefore, the order of
the kernel $(U_{\ell}/U_{\ell}^3)_{\b}$ of this homomorphism equals $\#
\mu(\ell)_3$. 

\vskip1mm
  
The short exact sequence $(4)$ in the proof of Proposition 0.12 of  [BP]
gives us the following exact sequence: $$1\rightarrow
(U_{\ell}/U_{\ell}^3)_{\b}\rightarrow
(\ell_3/{\ell^{\times}}^3)_{\b}\rightarrow (\cP\cap \cI^3)/\cP^3,$$
where $(\ell_3/{\ell^{\times}}^3)_{\b} = (\ell_3/{\ell^{\times}}^3)\cap
(\ell^{\times}/{\ell^{\times}}^3)_{\b}$, $\cP$ is the group of all
fractional 
principal ideals of $\ell$, and $\cI$ the group of all fractional ideals (we use multiplicative
notation for the group operation in both $\cI$ and $\cP$). Since the
order of the last group of the above exact sequence is $h_{\ell,3}$,
see $(5)$ in the proof of Proposition 0.12 of  [BP], we conclude that 
$$\#(\ell_3/{\ell^{\times}}^3)_{\b} \leqslant \#\mu(\ell)_3\cdot h_{\ell,3}.$$

Now we note that the order of the first term of the short exact sequence 
of Proposition 2.9 of [BP], for $G' =G$ and $S =V_{\infty}$, 
is $3^r /\#\mu(\ell)_3$.

Using the above observations, together with Proposition 2.9 and Lemma
5.4 
of [BP], and a
close 
look at the 
arguments in 5.3 and 5.5 of [BP] for $S=V_{\infty}$ and $G$ as above, we can derive the following upper bound: 
\begin{equation}
[\Gamma : \Lambda] \leqslant h_{\ell, 3}3^{r+\#\cT}.  
\end{equation}
This, together with (41) leads to the following bound:
\begin{equation}
(D_kD_{\ell})^{13} < \Big( \frac{(2\pi)^{42}}{4!5!7!8!11!}\Big)^dh_{\ell}.
\end{equation}

\ni{\bf 6.3.} Let $$A=\frac{(2\pi)^{42}}{4!5!7!8!11!}.$$  From bound (42), using (30), we obtain $$(D_kD_{\ell})^{13}< h_{\ell}A^d\leqslant A^d\frac{\delta(1+\delta)\Gamma(1+\delta)^d D_{\ell}^{(1+\delta)/2}\zeta_{\ell}(1+\delta)}{(R_{\ell}/w_{\ell})(2\pi)^{d(1+\delta)}}.$$ 
Hence, \begin{equation}
D_k^{13}D_{\ell}^{13-\frac{1+\delta}{2}}< A^d\frac{\delta(1+\delta)\Gamma(1+\delta)^d \zeta_{\ell}(1+\delta)}{(R_{\ell}/w_{\ell})(2\pi)^{d(1+\delta)}}.\end{equation}
As $D_k^2\leqslant D_{\ell}$, and $\zeta_{\ell}(1+\delta)\leqslant \zeta(1+\delta)^{2d}$, we conclude that $$D_k^{38-\delta}< A^d\frac{\delta(1+\delta)\Gamma(1+\delta)^d \zeta(1+\delta)^{2d}}{(R_{\ell}/w_{\ell})(2\pi)^{d(1+\delta)}}.$$  Therefore,
\begin{equation}
D_k^{1/d}<\Big[\{A\frac{\Gamma(1+\delta)\zeta(1+\delta)^2}{(2\pi)^{1+\delta}}
\}\cdot
\{\frac{\delta(1+\delta)}{R_{\ell}/w_{\ell}}\}^{1/d} \Big]^{1/(38-\delta)}.\end{equation}
Using the lower bound  
$R_\ell\geqslant 0.02 w_\ell\,e^{0.1d}$ due to R.\,Zimmert [Z], we obtain from this  the following:
\begin{equation}
D_k^{1/d} <f(d,\delta):=\Big[\{A\frac{\Gamma(1+\delta)\zeta(1+\delta)^2}{(2\pi)^{1+\delta}e^{0.1}}
\}\cdot
\{50\delta(1+\delta)\}^{1/d} \Big]^{1/(38-\delta)}.
\end{equation}

\ni From bound (43) we also obtain, 
\begin{equation}
D_\ell/D_k^2< \Big[ A
^dD_k^{-39}h_{\ell}\Big]^{1/13}.
\end{equation}

\ni Furthermore, using (44) and Zimmert's bound  $R_{\ell}\geqslant 0.02w_{\ell}e^{0.1d}$, we get the following:
%$$D_\ell/D_k^2\leqslant {\mathfrak p}_1(d,D_k,h_\ell):=\Big[ \Big(3A
%\Big)^d\frac{h_\ell}{D_k^{39}}\Big]^{1/13}$$

%and
\begin{equation}
D_\ell/D_k^2 < {\mathfrak p}(d,D_k,\delta)
:=\Big[\{A\frac{\Gamma(1+\delta)\zeta(1+\delta)^2}{(2\pi)^{1+\delta}e^{0.1}}
\}\cdot
\{\frac{50\delta(1+\delta)}{D_k^{38-\delta}}\}^{1/d} \Big]^{2d/(25-\delta)}.
\end{equation}
\bs

\ni {\bf 6.4.} For a fixed $\delta\geqslant 0.02$, $f(d,\delta)$ clearly decreases as $d$ increases. 
For $d\geqslant2,$ $D_k^{1/d} < f(d,2)\leqslant f(2,2)<2.3.$  We conclude now  from Proposition 2 that $d\leqslant 2,$ and for $d=2,$
$D_k\leqslant 5.$  Then $D_k=5.$ It
follows from bound (48) that 
$D_\ell/D_k^2< {\mathfrak p}(2,5,2)<2.$   Hence, $D_\ell/D_k^2=1$ and $D_\ell=25,$ which contradicts the
bound given by Proposition 2. We conclude that $d=1$, i.e., $k=\bQ.$
\vskip1mm

It is known, and follows, for example, from Proposition 7.1 of [PR], that a $\bQ$-group $G$ of type $^2E_6$, which at the unique real place of $\bQ$ is the outer form of rank $2$ (this is the form $^2E^{16'}_{6,2}$ which gives rise to a hermitian symmetric space), is isotropic over $\bQ$. This contradicts the fact that $G$ is anisotropic over $\bQ$ ({1.5}), and hence we conclude that groups of type $^2E_6$ do not give rise to arithmetic fake compact hermitian symmetric spaces. 

\vskip5mm

\begin{center}
{\bf 7. $G$ of type ${E}_7$ }  
\end{center}

\bs
\ni{\bf 7.1.}  In this section $G$ is assumed to be of type $E_7$. The dimension of $G$ is $133$, the exponents of its Weyl group are $1$, $5$, $7$, $9$, $11$, $13$ and $17$; and $s =2$. The dimension of the symmetric space $X$ of $\cG =\prod_{j=1}^rG(k_{v_j})$ is $ 27r$.  Let $\cT$ be as in 2.10. The bound (13) in this case gives us the following:
\begin{equation}
D_k^{133/2}< \frac{[\Gamma:\Lambda]}{2^{r+\#\cT}}\cdot\Big(\frac{(2\pi)^{70}}{5!7!9!11!13!17!}\Big)^d.
\end{equation}
\vskip1mm

The center $C$ of $G$ is $k$-isomorphic to $\mu_2$. 
The Galois cohomology group $H^1(k, C)$ is isomorphic to $k^{\times}/{k^{\times}}^{2}$. The order of the first term of the short exact sequence of Proposition 2.9 of [BP], for $G' =G$ and $S = V_{\infty}$, is $2^{r-1}$. From the proof of Proposition 0.12 of [BP], we easily conclude that $\#k_2/{k^{\times}}^2 \leqslant h_{k,2}2^d$. We can adapt the argument used to prove Proposition 5.1 in [BP], and the argument in 5.5, of [BP], for $S=V_{\infty}$ and $G' =G$, to derive the following bound:
\begin{equation}
[\Gamma : \Lambda]\leqslant h_{k,2}2^{d+r-1+\#T}.
\end{equation}  
Combining (49) and (50) we obtain the following bound:
\begin{equation}
D_k^{133/2} < 2^{d-1}\Big(\frac{(2\pi)^{70}}{5!7!9!11!13!17!}\Big)^dh_{k,2}.
\end{equation}
\vskip1mm

\ni{\bf 7.2.} Let $$B = \frac{(2\pi)^{70}}{5!7!9!11!13!17!}.$$  From (51) we obtain the following:
$$D_k^{1/d} <\Big[ 2B
(h_{k,2}/2)^{1/d}\Big]^{2/133}.$$
Using the Brauer-Siegel bound (16) for totally real number fields, and the obvious bound $\zeta_k(1+\delta)\leqslant \zeta(1+\delta)^d$, we obtain  
\begin{equation}
D_k^{1/d}<\Big[\{B\frac{\Gamma((1+\delta)/2)\zeta(1+\delta)}{\pi^{(1+\delta)/2}}
\}\cdot
\{\frac{\delta(1+\delta)}{R_{k}}\}^{1/d} \Big]^{2/(132-\delta)}.
\end{equation}
Now using the lower bound  
$R_k\geqslant 0.04\,e^{0.46d}$ due to R.\,Zimmert [Z] again, we get
\begin{equation}
D_k^{1/d} <\phi (d,\delta):=\Big[\{B\frac{\Gamma((1+\delta)/2)\zeta(1+\delta)}{\pi^{(1+\delta)/2}e^{0.46}}
\}\cdot
\{25\delta(1+\delta)\}^{1/d} \Big]^{2/(132-\delta)}.
\end{equation}

\bs
\ni{\bf 7.3.}  For a fixed $\delta\geqslant 0.04$, $\phi (d,\delta)$ clearly decreases as $d$ increases.  
By a direct computation we see that $\phi (2,4)< 2$, and hence  for all totally real  number field $k$
of degree $d\geqslant 2,$
$$D_k^{1/d} < \phi (d,4)\leqslant \phi (2,4)<2.$$
From this bound and Proposition 2 we conclude that $d$ can only be $1$, i.e., $k =\bQ$. But then $r = 1$ and the complex dimension of the associated symmetric space $X$ is $27$. Then the Euler-Poincar\'e characteristic of any quotient of $X$ by a cocompact torsion-free discrete subgroup of $\oG$ is negative (1.3), and hence it cannot be a fake compact hermitian symmetric space.  Another way to eliminate this case is to observe that an absolutely simple $\bQ$-group of type $E_7$ is isotropic if it is isotropic 
over $\bR$ (this result follows, for example, from Proposition 7.1 of [PR]). 

\vskip6mm

\begin{center}
{\bf 8. $G$ of  type $^2A_n$ with $n$ odd}  
\end{center}
\vskip4mm

\ni{\bf 8.1.} We shall assume from now on that $G$ is an absolutely simple simply connected $k$-group of type $^2A_n$ with $n>1$ odd. We retain the notation introduced in \S\S 1,\,2. 
 In particular,  $\ell$ is the totally complex quadratic extension of $k$ over which $G$ is an inner form, $d = [k:\bQ]$, $s = n+1$; $\Gamma$, $\Lambda$,  for $v\in V_f$, the parahoric subgroups $P_v$ of $G(k_v)$ are as in 2.1, and $\cT$ is as in 2.10. We recall that  for every nonarchimedean $v\notin \cT$, $\Xi_{\Theta_v}$ is trivial. For all $v\in \cT$,  $\# \Xi_{\Theta_v}| (n+1)$ and $e(P_v)>e'(P_v)>n+1$. We also recall from 2.1 that $\mu(\cG/\Gamma)$ is a submultiple of $1/(n+1)^r$, hence, $(n+1)^r\mu(\cG/\Gamma)\leqslant 1$.
\vskip1mm

The center $C$ of $G$ is the kernel of the norm map $N_{\ell/k}: R_{\ell/k}(\mu_{n+1})\to \mu_{n+1}$. Therefore, we get the following exact sequence:$$(*)\ \ \ \ \ \ \ \ \ \ \ 1\to \mu_{n+1}(k)/N_{\ell/k}(\mu_{n+1}(\ell))\to H^1(k, C) \to (\ell^{\times}/{\ell^{\times}}^{n+1})_{\bullet}\to 1,$$ where $(\ell^{\times}/{\ell^{\times}}^{n+1})_{\bullet}$ is the kernel of the homomorphism $\ell^{\times}/{\ell^{\times}}^{n+1}\to k^{\times}/{k^{\times}}^{n+1} $ induced by the norm map $N_{\ell/k}:\ell^{\times}\to k^{\times}$.  By Dirichlet's unit theorem, $U_k\cong \{\pm 1\}\times \bZ^{d-1}$ and $U_{\ell}\cong \mu(\ell)\times \bZ^{d-1}$, and hence, $U_k/U_k^{n+1}\cong \{\pm 1\}\times (\bZ/(n+1)\bZ)^{d-1}$ and $U_{\ell}/U_{\ell}^{n+1}\cong \mu_{n+1}(\ell)\times (\bZ/(n+1)\bZ)^{d-1}$.  Since $N_{\ell/k}(U_{\ell}) \supset N_{\ell/k}(U_k) = U_k^2$, the image of the homomorphism $U_{\ell}/U_{\ell}^{n+1}\to U_k/U_k^{n+1}$ induced by the norm map $N_{\ell/k}$ contains $U_k^2/U_k^{n+1}\, (\cong (2\bZ/(n+1)\bZ)^{d-1})$, and hence the kernel $(U_{\ell}/U_{\ell}^{n+1})_{\bullet}$ of this homomorphism is of order at most $\#\mu_{n+1}(\ell)\cdot 2^{d-1}$.  The short exact sequence (4) in the proof of Proposition 0.12 of [BP] gives us the following exact sequence:$$1\to (U_{\ell}/U_{\ell}^{n+1})_{\bullet}\to (\ell_{n+1}/{\ell^{\times}}^{n+1})_{\bullet}\to (\cP\cap \cI^{n+1})/\cP^{n+1},$$ where $\ell_{n+1}$ is the subgroup of $\ell^{\times}$ consisting of all $x$ such that for every normalized nonarchimedean valuation $v$ of $\ell$, $v(x)\in (n+1)\bZ$, $(\ell_{n+1}/{\ell^{\times}}^{n+1})_{\bullet} = (\ell_{n+1}/{\ell^{\times}}^{n+1})\cap (\ell^{\times}/{\ell^{\times}}^{n+1})_{\bullet}$, $\cP$ is the group of all fractional principal ideals of $\ell$, and $\cI$ the group of all fractional ideals (we use multiplicative notation for the group operation in both $\cI$ and $\cP$). Since the order of the last group of the above exact sequence is $h_{\ell, n+1}$, see (5) in the proof of Proposition 0.12 of [BP], we conclude that $\#(\ell_{n+1}/{\ell^{\times}}^{n+1})_{\bullet}\leqslant \#\mu_{n+1}(\ell)\cdot 2^{d-1}h_{\ell, n+1}$. 
\vskip1mm

Let $c$ be the order of the kernel of the norm map  $N_{\ell/k}: \mu_{n+1}(\ell)\to \mu_{n+1}(k)=\{\pm 1\}$. Then the order of the first term of $(*)$ is $2c/\#\mu_{n+1}(\ell)$, whereas the order of the first term of the short exact sequence of Proposition 2.9 of [BP],  for $G' = G$ and $S =V_{\infty}$, is $(n+1)^r/c$. Now from Lemma 5.4 of [BP] and the arguments given in 5.3 and 5.5 of that paper  we obtain the following upper bound (note that we need to replace ``$n$" in 5.3 and 5.5 of [BP] with ``$n+1$" since the group $G$ in this and the next section  is of type $^2A_n$): 
\begin{equation}
[\Gamma:\Lambda]\leqslant h_{\ell,n+1}2^d(n+1)^{r+\#\cT}.
\end{equation}

For the group $G$ under consideration here, $\dim\,G = n^2+2n$, the exponent $m_j = j$ and $\mathfrak{s} = (n-1)(n+2)/2$, so the volume formula (3) gives us the following: 
\begin{equation}
\mu(\cG /\Lambda )= D_k^{\frac{1}{2}{(n^2+2n)}}(D_{\ell}/D_k^2)^{\frac{1}{4}{(n-1)(n+2)}}\Big(\prod_{j=1}^{n}\frac{j!}{(2\pi)^{j+1}}\Big)^{d}\prod_{v\in V_f}e(P_v).
\end{equation}
For $v\in \cT$, as $e(P_v)>(n+1)$, and moreover for all $v\in V_f$, $e(P_v)>1$, using (54) and (55) we find that 
\begin{equation}
1\geqslant (n+1)^r\mu(\cG /\Gamma )>D_k^{\frac{1}{2}{(n^2+2n)}}(D_{\ell}/D_k^{2})^{\frac{1}{4}{(n-1)(n+2)}}\Big(\prod_{j=1}^{n}
\frac{j!}{(2\pi)^{j+1}}\Big)^{d}\frac{1}{2^dh_{\ell, n+1}},
\end{equation}
As $D_k^2 |D_{\ell}$, from (56) we obtain  the following bound  for $D_k$:
\begin{equation}
D_k^{1/d} <
f_1(n,d,h_{\ell, n+1}):=\Big[ \{2
\prod_{j=1}^{n}\frac{(2\pi)^{j+1}}{j!}\}^d\cdot h_{\ell, n+1}
\Big]^{\frac{2}{d(n^2+2n)}}.
\end{equation}
Since $\zeta_{\ell}(1+\delta)\leqslant \zeta(1+\delta)^{2d}$, for $\delta>0$, we obtain the following bound from (57)  and (30)

\begin{eqnarray}
D_k^{1/d}&<&f_2(n,d,R_{\ell}/w_{\ell},\delta)\\
\nonumber&:=&\Big[\{2\frac{\Gamma(1+\delta)\zeta(1+\delta)^2}{(2\pi)^{1+\delta}}
\prod_{j=1}^{n}\frac{(2\pi)^{j+1}}{j!}\}\cdot
\{\frac{\delta(1+\delta)}{({R_{\ell}}/{w_{\ell}})}\}^{1/d} \Big]^{\frac{2}{(n^2+2n-2\delta-2)}}.
\end{eqnarray}

Using the lower bound  
$R_\ell\geqslant 0.02 w_\ell\,e^{0.1d}$ due to Zimmert, we obtain the following from (30)
\begin{equation}
\frac{1}{h_{\ell,n+1}}\geqslant\frac1{h_\ell}\geqslant\frac{0.02}{\delta(1+\delta)}
\Big( \frac{(2\pi)^{1+\delta} e^{0.1}}{\Gamma(1+\delta)}\Big)^d
\frac1{D_\ell^{(1+\delta)/2}\zeta_\ell(1+\delta)}.
\end{equation}

\ni Since $D_\ell\geqslant D_k^2,$ and $\zeta_\ell(1+\delta)\leqslant \zeta(1+\delta)^{2d},$ we get the following for all $\delta$ in the interval $ [0.02,6.5]$.  
\begin{equation}
D_k^{1/d} <f_3(n,d,\delta):=\Big[ \{2 \frac{\Gamma(1+\delta)\zeta(1+\delta)^2}{(2\pi)^{1+\delta}
e^{0.1}}
\prod_{j=1}^{n}\frac{(2\pi)^{j+1}}{j!}\}\cdot
\{50\delta (1+\delta) \}^{1/d} \Big]^{\frac{2}{(n^2+2n-2\delta-2)}}.
\end{equation}

Now the following three bounds for the relative discriminant $D_{\ell}/D_k^2$ are obtained from (56), (59) and (30).

\begin{equation}
D_{\ell}/D_k^2< {\mathfrak p}_1(n,d,D_k,h_{\ell,n+1}):=\Big[h_{\ell,n+1}\cdot \{ 
2\prod_{j=1}^{n}\frac{(2\pi)^{j+1}}{j!}\}^d D_k^{-(n^2+2n)/2}\Big]^{\frac{4}{(n-1)(n+2)}}.
\end{equation}

\begin{eqnarray} 
\ \ \ \ \ \ \ \ \ \ \ \ D_{\ell}/D_k^2&<& {\mathfrak p}_2(n,d,D_k,R_{\ell}/w_{\ell},\delta)\\ 
\nonumber&:=& \Big[\frac{\delta(1+\delta)}{(R_\ell/w_\ell)D_k^{(n^2+2n-2\delta-2)/2}}
 \{ \frac{2\Gamma(1+\delta)\zeta(1+\delta)^2}{(2\pi)^{1+\delta}}
\prod_{j=1}^{n}\frac{(2\pi)^{j+1}}{j!}\}^d \Big]^{\frac{4}{(n^2+n-2\delta-4)}}.
\end{eqnarray}

\begin{eqnarray}
D_\ell/D_k^2 &<&{\mathfrak p}_3(n,d,D_k,\delta )\\
\nonumber&:=&\Big[\frac{50\delta(1+\delta)}{D_k^{(n^2+2n-2\delta -2)/2}}\cdot \{ 
\frac{2\Gamma(1+\delta)\zeta(1+\delta)^2}{(2\pi)^{1+\delta }e^{0.1}}
\prod_{j=1}^{n}\frac{(2\pi)^{j+1}}{j!}\}^d \Big]^{\frac{4}{(n^2+ n-2\delta-4)}}.\end{eqnarray}

We also get the following bound for $D_\ell$ from  (56).
\begin{equation}
D_{\ell}^{1/2d}< {\mathfrak q}_1(n,d,D_k,h_{\ell,n+1}):=\Big[\frac{h_{\ell,n+1}}{D_k^{(n+2)/2}} \cdot \{ 
2\prod_{j=1}^{n}\frac{(2\pi)^{j+1}}{j!}\}^d \Big]^{\frac2{d(n-1)(n+2)}},
\end{equation}
which in turn provides the following bound using (30) and (59)
\begin{eqnarray} 
D_{\ell}^{1/2d}&<& {\mathfrak q}_2(n,d,D_k,R_{\ell}/w_{\ell},\delta) \\
\nonumber&:=& \Big[\frac{\delta(1+\delta)}{(R_\ell/w_\ell)D_k^{(n+2)/2}}
\cdot \{ \frac{2\Gamma(1+\delta)\zeta(1+\delta)^2}{(2\pi)^{1+\delta}}\cdot
\prod_{j=1}^{n}\frac{(2\pi)^{j+1}}{j!}\}^d \Big]^{\frac2{d(n^2+n-2\delta-4)}}.\end{eqnarray}

\vskip4mm

We now state the following simple lemma.

\begin{lemm}
Let $\delta\in [0.02,6.5].$  For fixed values of $n$ and $\delta,$ 
$f_3(n,d,\delta)$ decreases as $d$ increases.  Furthermore, for fixed values of $d$ and $\delta,$ if 
$n \geqslant 7,$  then $f_3(n+1,d,\delta)<\max(1, f_3(n,d,\delta)).$ 
\end{lemm}

\bs
\ni{\bf 8.2.}  
Let us begin determination of the totally real number field $k.$
Let $f_3(n,d,\delta)$ be as in $(60)$. By a direct computation we obtain the following  
upper bound for the value of $f_3(n,2,3)$ for small $n$.

$$\begin{array}{ccc}
n:&13&11\\
f_3(n,2,3)<&2.1&2.4.
\end{array}$$
Hence for $n\geqslant 13$ and $d\geqslant 2$,  $f_3(n,d,3)\leqslant f_3(13,2,3)<2.1,$ 
which in view of Proposition 2 implies that $k =\bQ$. 

\vskip1mm
 
\ni{\bf 8.3.} Now we will determine the degrees $d$ of 
possible $k$ for $n\leqslant 11$ using $(60).$  We get the following table by evaluating $f_3(n,d,\delta)$, with $n$ given in the first column, $d$ given in the second column, and $\delta$ given in the third column

$$\begin{array}{cccc}
n&d&\delta&f_3(n,d,\delta)<\\
11&3&2&2.4\\
9&3&1.9&2.9\\
7&3&1.6&3.63\\
5&4&1.3&5.12\\
3&7&1&10.
\end{array}$$
Taking into account the upper bound in the last column of the above table, Proposition 2 implies the following bound for
$d$ for each odd integer $n$ between $3$ and $11.$
\vskip1.5mm

$$\begin{array}{cccccc}
n:&11&9&7&5&3\\
d\leqslant&2&2&2&3&6.
\end{array}$$

\bs
\ni{\bf 8.4.}  We will now narrow down the possibilities for 
$d$ further.  We begin with larger values
of $n.$  
\vskip1mm

For $n=11,9$ and $7$, we know that $d\leqslant 2.$ 
\vskip1mm
 
For $n=11$ and $d=2,$  
$D_k^{1/2}\leqslant f_3(11,2,2)\leqslant 2.5$, so $D_k=5.$ Then 
$D_\ell/D_k^2\leqslant\lfloor {\mathfrak p}_3(11,2,5,2)\lr$ $=1.$  Hence $D_\ell=25,$ 
but there is no such $\ell$. This implies that if  $n =11$, then $k =\bQ$.
\vskip1mm

For $n=9$ and $d=2,$ 
$D_k^{1/2}\leqslant f_3(9,2,2))<3.$
Hence, $D_k=5$ or $8.$  As 
$\lfloor{\mathfrak p}_3(9,2,5,2)\lr=3$ and $\lfloor{\mathfrak p}_3(9,2,8,1.6)\lr=1.$
So $D_\ell\leqslant75,$ but there is no $\ell$ of degree $4$ for which this bound holds, 
and we conclude that if $n=9$, then again $k =\bQ$.
\vskip2mm

\ni{\bf 8.5.} We shall now consider the case $n=7$ and $d=2$. As 
$D_k^{1/2}\leqslant f_3(7,2,1)<3.8$, 
$D_k=5,8,12$ or $13$. 
Computations give that $\lfloor{\mathfrak p}_3(7,2,5,1.3)\lr=11,$
 $\lfloor{\mathfrak p}_3(7,2,8,1.3)\lr=3,$ $\lfloor{\mathfrak p}_3(7,2,12,1.3)\lr=1$
 and $\lfloor{\mathfrak p}_3(7,2,13,1)\lr=1$
Hence $D_\ell$ is bounded from above by $\max(5^2\cdot 11, 8^2\cdot 3,12^2, 13^2).$  From the list of number fields given in [1], we conclude that
the class number of all these totally complex quartic $\ell$ is $1$.  Hence the pairs $(k,\ell)$ belong to the list of 
[PY1], 8.2 (see also  [PY1], 7.10). 
Also, the bound for 
the relative discriminant $D_{\ell}/D_k^2$ can  be improved
to
$\lfloor{\mathfrak p}_1(7,2,5,1)\lr=8,$ and
 $\lfloor{\mathfrak p}_1(7,2,8,1)\lr=2$ in the first two cases.  Now checking against the list
of [PY1], 8.2,  we conclude that the following are the only possible pairs $(k,\ell)$.
$$\cC_1,\: \cC_{11}.$$ 
We eliminate these pairs by computing $\cR$ and then using Proposition 1. The values of $\zeta_k$ and $\zeta_{\ell | k}$ required for the computation of $\cR$ are given below.
 $$\begin{array}{ccccccccl}
(k,\ell)&\zeta_k(-1)&\zeta_{\ell|k}(-2)&\zeta_k(-3)&\zeta_{\ell|k}(-4)&\zeta_k(-5)&\zeta_{\ell|k}(-6)&\zeta_k(-7)\\
\cC_1&1/30&4/5&1/60&1172/25&67/630&84676/5&361/120\\
\cC_{11}&1/6&1/9&23/60&5/3&1681/126&427/3&257543/120.\\
\end{array}$$
\vskip1mm

\ni{\bf 8.6.} Consider now $n=5.$  We know that $d\leqslant 3.$ Assume, if possible,  that $d =3$. 
As $D_k^{1/3}< f_3(5,3,1)<5.3$,  we see from the table of totally real cubics given in [C]
that $D_k$  is either $49$ or $81$. On the other hand, 
$D_\ell/D_k^2\leqslant\lfloor{\mathfrak p}_3(5,3,49,1)\lr=16$
and $D_\ell/D_k^2\leqslant\lfloor{\mathfrak p}_3(5,3,81,1)\lr=4$ for the two cases respectively.  So  
$D_\ell\leqslant \max(49^2\cdot 16, 81^2\cdot 4)=49^2\cdot 16.$  The 
class number of all totally complex sextic fields  $\ell$ with $D_\ell\leqslant 49^2\cdot 16$ is $1$.  Now the bound 
for the relative discriminant $D_{\ell}/D_k^2$ can be improved to  
$\lfloor{\mathfrak p}_1(5,3,49,1)\lr=6,$ and
 $\lfloor{\mathfrak p}_1(5,3,81,1)\lr=1$ in the two cases.  Among the pairs $(k,\ell)$ listed in 
 [PY1], 8.2,  none satisfy these  
conditions.  Hence, $d<3.$
\vskip1mm

Assume now $d=2.$  As $D_k^{1/2}< f_3(5,2,1)<5.54$, $D_k\leqslant30.$ From Friedman [F] we know that $R_\ell/w_\ell\geqslant 1/8$
except when $D_\ell=117, 125$ and $144$.  Therefore, apart from the 
three exceptional cases, we conclude that
$D_\ell/D_k^2\leqslant \lfloor{\mathfrak p}_2(5,2,5,1/8,1)\lr=82.$  Since the discriminant  in each of the
three exceptional cases is smaller than $82\cdot 5^2,$ we conclude that
the bound $D_\ell/D_k^2\leqslant 82$ always holds.  
So 
$D_\ell\leqslant30^2\cdot 82=73800.$  From the list in [1] of totally complex quartics 
$\ell$ with $D_\ell\leqslant 73800,$ we see  that $h_\ell\leqslant15$ and 
hence $h_{\ell,6}\leqslant 12.$
Then  $D_k^{1/2}< f_1(5,2,12)<5.1$, and so $D_k\leqslant 24$. 
We know that either $D_k=5$ or $D_k\geqslant 8.$  In  the latter case, as $\ll{\mathfrak p}_1(5,2,8,12)\lr=17.$ 
Thus $D_\ell\leqslant\max(5^2\cdot82, 24^2\cdot 17)=9792.$  By checking the list of totally complex quartic number fields in [1] again, we conclude that
$h_\ell\leqslant 5$ and hence $h_{\ell,6}\leqslant 4.$
Then  $D_k^{1/2}< f_1(5,2,4)<4.87$, so $D_k\leqslant 21.$  We now compute 
$\lfloor{\mathfrak p}_1(5,2,D_k,4)\lr$ for $5\leqslant D_k\leqslant 21$ to get the following  bound for $D_{\ell}$:
$$\begin{array}{ccccccc}
D_k:&5&8&12&13&17&21\\
\lfloor{\mathfrak p}_1(5,2,D_k,4)\lr:&48&15&5&4&2&1\\
D_\ell\leqslant&1200&960&720&676&578&441.
\end{array}$$
The list of  number fields satisfying the above constraint was provided by Malle using the tables in [1].  
It turns out that all the number fields involved
are listed in the tables in [PY1], 8.2.  Moreover, all have class number $1.$  It follows that $D_\ell$ is bounded by 
 $\lfloor{\mathfrak p}_1(5,2,D_k,1)\lr$.
 $$\begin{array}{ccccccc}
D_k:&5&8&12&13&17&21\\
\lfloor{\mathfrak p}_1(5,2,D_k,1)\lr:&40&12&4&3&1&1\\
D_\ell\leqslant&1000&768&576&507&289&441.
\end{array}$$ From  the table in [PY1], 8.2, we conclude that the following are the
ony possibilities for the pair $(k,\ell)$.
$$ \cC_1,\:\cC_2,\:\cC_3,\: \cC_8,\:\cC_{9},\:\cC_{11},\:\cC_{17}.$$ 
We eliminate each of the above pairs by computing $\cR$, using the following values of $\zeta_k$ and $\zeta_{\ell |k}$, and then use Proposition 1. 

$$\begin{array}{ccccccccl}
(k,\ell)&\zeta_k(-1)&\zeta_{\ell|k}(-2)&\zeta_k(-3)&\zeta_{\ell|k}(-4)&\zeta_k(-5)\\
\cC_1&1/30&4/5&1/{60}&1172/{25}&{67}/{630}\\
\cC_2&1/{30}&{32}/9&1/{60}&{1984}/3&{67}/{630}\\
\cC_3&1/{30}&15&1/{60}&8805&{67}/{630}\\
\cC_8&1/{12}&3/2&{11}/{120}&{285}/2&{361}/{252}\\
\cC_9&1/{12}&{92}/9&{11}/{120}&{15940}/3&{361}/{252}\\
\cC_{11}&1/6&1/9&{23}/{60}&5/3&{1681}/{126}\\
\cC_{17}&1/3&{32}/{63}&{77}/{30}&{64}/3&{17971}/{63}.\
\end{array}
$$

\vskip3mm
 
\ni{\bf 8.7.} The case $n=3$ requires more detailed considerations.
\vskip1mm

\ni$\bullet$ Again we are considering totally real $k$ with $d>1$ in this section. 
We know from 8.3 that 
$d\leqslant 6$.  Consider first $d=6.$  
Then $D_k\geqslant 300125$ (see \S 3). 
Hence,  $D_\ell^{1/12}\leqslant {\mathfrak q}_2(3,6,300125,1/8,1)<12.$
 According to Table IV of [Mart],
$M_c(32)> 12$, so considering the Hilbert class field of $\ell$ which is an extension of degree $h_{\ell}$ of $\ell$, 
 we infer that 
$h_\ell\leqslant\ll31/12\lr=2.$
Hence $h_{\ell,n+1}\leqslant 2.$ Now applying bound  (64) we obtain  
$D_k^{1/6}\leqslant D_\ell^{1/12}\leqslant {\mathfrak q}_1(3,6,300125,2)<7,$ 
which contradicts Proposition 2.

\vskip2mm
\ni$\bullet$ Consider now $d=5.$  In this case $D_k\geqslant14641.$ 
Hence, $D_\ell^{1/10}\leqslant {\mathfrak q}_2(3,5,14641,1/8,1)$ $<14.$
According to Table IV of [Mart],
$M_c(52)> 14.1.$  Using again the Hilbert class field of $\ell$ we conclude that  
$h_\ell\leqslant\ll51/10\lr=5$, and hence $h_{\ell,4}\leqslant 4.$
Then  $D_\ell^{1/10}\leqslant {\mathfrak q}_1(3,5,14641,4)<7.74$
So $D_\ell<7.74^{10}<7.72\times 10^8.$  On the other hand, Schehrazad Selmane [Sel] has shown that  the totally complex number field of degree $10$, containing a totally real quintic field, with smallest absolute discriminant is the cyclotomic field 
$\bQ(\zeta_{11})$ generated by a primitive $11$-th root $\zeta_{11}$ of unity. This field has absolute discriminant $11^9$. Since 
$11^9> 7.72\times 10^8$, we conclude that $d\ne 5$.   
\vskip2mm

\ni$\bullet$ Consider now $d=4.$  In this case $D_k\geqslant725.$
Hence $D_\ell^{1/8}\leqslant {\mathfrak q}_2(3,4,725,1/8,0.86)<17.43.$
According to Table IV of [Mart],
$M_c(140)> 17.49$,  so considering the Hilbert class field of $\ell$ we find that 
$h_\ell\leqslant\ll139/8\lr=17.$
So  $h_{\ell, 4}\leqslant 16$ and 
then $D_\ell^{1/8}\leqslant {\mathfrak q}_1(3,4,725,16)<9.7.$
According to Table IV of [Mart],
$M_c(20)>9.8$ which by considering the Hilbert class field of $\ell$ implies that  
$h_\ell\leqslant\ll19/8\lr=2$ and $h_{\ell,4}\leqslant 2.$  It follows that
$D_\ell^{1/8}\leqslant {\mathfrak q}_1(3,4,725,2)<8.7.$
According to Table IV of [Mart],
$M_c(16)>8.7.$  Hence, again by considering the Hilbert class field of $\ell$ implies that 
$h_\ell\leqslant\ll15/8\lr=1.$  But then  $D_\ell^{1/8}\leqslant {\mathfrak q}_1(3,4,725,1)<8.386.$
So $D_\ell\leqslant \ll8.386^8\lr<2.45\times 10^7.$
Also, $D_k^{1/4}<f_1(3,4,1)<7.146.$  Hence, $D_k\leqslant 2607.$
We also know that 
$D_\ell/D_k^2\leqslant \ll{\mathfrak p}_1(3,4,725,1)\lr=46.$
Any such pair $(k,\ell)$ lies in the list of pairs tabulated in [PY1], 8.2. We find that the possible pairs are $\cC_{34}-\cC_{37}$ in the notation used in [PY1], 8.2.
Again, we eliminate each of the pairs by computing $\cR$ using the following zeta values and applying Proposition 1.
$$\begin{array}{cccc}
(k,\ell)&\zeta_k(-1)&\zeta_{\ell|k}(-2)&\zeta_k(-3)\\
\cC_{34}&4/15&128/45&2522/15\\
\cC_{35}&2/3&12&3793/3\\
\cC_{36}&5/6&411&87439/60\\
\cC_{37}&1&46/3&22011/10.\end{array}
$$
\ni$\bullet$ Consider now $d=3.$  The three smallest absolute discrimants of 
totally real cubic fields are $49, 81$ and $148.$  Let us consider first 
the totally real cubic fields $k$ with $D_k\geqslant 148.$   
Note that $D_\ell^{1/6 }<{\mathfrak q}_2(3,3,148,1/8,0.7)<18.1$,
since $R_\ell/w_\ell\geqslant 1/8$ except for the six sextics whose discriminants are  listed in
[PY1],\:7.3. The  root discriminant of these six sextics clearly satisfy the above bound.
We see from Table IV in [Mart] that $M_c(180)>18.1.$  Hence, considering the Hilbert class field of $\ell$, 
we conclude that $h_\ell\leqslant \ll179/6\lr=29$ and so $h_{\ell,4}\leqslant 16.$
Then $D_\ell^{1/6}<{\mathfrak q}_1(3,3,148,16)<10.$
Again from Table IV of [Mart] we find that $M_c(22)>10.25$, and as before  
considering the Hilbert class field of $\ell$, we conclude that $h_\ell\leqslant \ll21/6\lr=3.$ So $h_{\ell,4}\leqslant 2.$
Then $D_\ell^{1/6}< {\mathfrak q}_1(3,3,148,2)<8.7.$ 
Furthermore, 
$D_k^{1/3}<f_1(3,3,2)<7.37.$
Hence $D_\ell\leqslant \ll8.7^6\lr=433626$ and $148\leqslant D_k\leqslant \ll7.37^3\lr=400$.  {There are only three
pair of number fields $(k,\ell)$ satisfying the above bounds and $h_\ell=1$ for each of the $\ell$ occurring in these three pairs from which we conclude that
$D_\ell^{1/6}< {\mathfrak q}_1(3,3,148,1)<8.31.$   Hence, $D_\ell<329311.$  But there are no pairs $(k,\ell)$ of totally real cubic $k$, and totally complex quadratic extension $\ell$ of $k$ with $148\leqslant D_k\leqslant 400$ and $D_{\ell}<329311$. }

\vskip1mm

We will consider now the unique totally real cubic field $k_1$ with $D_{k_1}=81.$  Note that  $k_1 =\bQ[x]/ (x^3-3x-1)$, the regulator 
$R_{k_1}\geqslant 0.849$ according to [C].  Now by listing $m$ such that  the value $\phi(m)$ of the Euler function $\phi$ is a divisor of 6,  we  see that unless $\ell$ is $\bQ(\zeta_{18})$, $w_{\ell} = 2,\, 4$ or $6$ (note that $\bQ(\zeta_{14})$ does not contain $k_1$).  As $\ell$ is a CM field which is a quadratic extension of $k_1$, $R_{\ell} = 2^{2}R_{k_1}/Q$, where $Q=1$ or $2$ (cf.\,[W]), hence unless $\ell$ is $\bQ(\zeta_{18})$,  
$R_\ell/w_\ell\geqslant 2\times0.849/6=0.283$. 
So either $\ell = \bQ(\zeta_{18})$ or 
$D_\ell^{1/6}< {\mathfrak q}_2(3,3,81,0.283,0.66)<19.4.$ From Table IV of [Mart] we find that  $M_c(340)>19.4.$  Hence, by considering the Hilbert class field of $\ell$, we successively get the following improved bounds for $h_{\ell, 4}$:  $h_\ell\leqslant \ll339/6\lr=56$, so 
$h_{\ell,4}\leqslant 32.$
 Therefore, $D_\ell^{1/6}<{\mathfrak q}_1(3,3,81,32)<11.6.$
Again in Table IV of [Mart] we see that  $M_c(30)>11.6.$  So $h_\ell\leqslant \ll29/6\lr=4$, and
$D_\ell^{1/6}<{\mathfrak q}_1(3,3,81,4)<10.1.$
But according to Table IV of [Mart], $M_c(22)>10.2.$  Hence, $h_\ell\leqslant \ll21/6\lr=3$, and $h_{\ell,4}\leqslant 2.$
It follows that $D_\ell/D_{k_1}^2\leqslant \ll{\mathfrak p}_1(3,3,81,2)\lr= 120,$ from which
we conclude that $D_\ell\leqslant 120\cdot 81^2=787320.$
Malle provided us a complete list of totally complex quadratic extensions $\ell$ of the above $k_1$ with $D_{\ell}\leqslant 787320$. This list consists of three fields whose absolute discriminant, defining monic polynomial and the value of $\zeta_{\ell |k_1}(-2)$  are given below. 
$$\begin{array}{cllll}
D_\ell&\ell&\zeta_{\ell |k_1}(-2)\\
19683&x^6-x^3+1&-104/27\\
419904&x^6+6x^4+9x^2+1&-7826\\
465831&x^6-3x^5+9x^4-13x^3+15x^2-9x+3&-10944.
\end{array}$$
The first of these fields is $\bQ(\zeta_{18})$. We shall denote the three pairs $(k_1,\ell)$ with $\ell$ from this list, and  $k_1 = \bQ[x]/(x^3-3x-1)$,  $\cE_1$, $\cE_2$ and $\cE_3$ respectively.
\vskip2mm

Let us now consider the unique totally real cubic field $k_2$ with $D_{k_2}=49$. Note that $k_2 =\bQ[x]/(x^3-x^2-2x+1)$, and from [C] we find that  its regulator  is larger than $0.525.$
Hence, as for $k_1$, we see that except for the cyclotomic field $\bQ(\zeta_{14})$ which has class number 1 and discriminant $7^5 = 16807$, $w_\ell =2$, $4$ or $6$, and for the noncyclotomic $\ell$, 
$R_{\ell}/w_{\ell}\leqslant 2\times 0.525/6=0.175.$ Therefore, $D_\ell/D_{k_2}^2\leqslant \ll{\mathfrak p}_2(3,3,49,0.175,0.64)\lr
=62697,$ and hence $D_{\ell}\leqslant 62697\cdot 49^2$.
Malle provided the authors a list of totally complex quadratic extensions $\ell$ of $k_2$ for which this bound holds. For every $\ell$ in this list, $h_\ell\leqslant 30$, and 
hence, $h_{\ell,4}\leqslant 16.$  Therefore, $D_\ell^{1/6}<{\mathfrak q}_1(3,3,49,16)<12.02.$ From Table IV in [Mart] we see that  $M_c(34)\geqslant 12.4.$  We conclude, as before, by considering the Hilbert class field of $\ell$, that
$h_\ell\leqslant \ll33/6\lr=5$,  so $h_{\ell,4}\leqslant 4.$  Then 
$D_\ell^{1/6}< {\mathfrak q}_1(3,3,49,4)<10.96.$  Hence $D_\ell <1.74\times 10^6.$
{From the list provided by Malle, we see  that there are eleven candidates for $\ell$, each with class number $h_\ell\leqslant 2.$  Therefore, $D_\ell<\ll{\mathfrak p}_1(3,3,49,2)\lr\times 49^2=1306144.$  For all totally complex quadratic extension $\ell$ of $k_2$  satisfying this bound, $h_\ell=1$ and
hence we conclude that $D_\ell<\ll{\mathfrak p}_1(3,3,49,1)\lr \times 49^2=991613.$  From the list provided by Malle, we see
that the possible $\ell$ are:}
$$\begin{array}{cllll}
D_\ell&\ell&\zeta_{\ell |k_2}(-2)\\
16807&x^6-x^5+x^4-x^3+x^2-x+1& -64/7\\
64827&x^6-x^5+3x^4+5x^2-2x+1&-2408/9\\
153664&x^6+5x^4+6x^2+1& -2306\\
400967&x^6-2x^5+5x^4-7x^3+10x^2-8x+8&-25536\\
573839&x^6-x^5+4x^4-3x^3+8x^2-4x+8&-62208\\
602651&x^6-3x^5+10x^4-15x^3+21x^2-14x+7&-70392\\
909979&x^6-2x^5+7x^4-12x^3+21x^2-15x+13&-196216.\\
\end{array}$$
We shall denote the seven pairs $(k_2,\ell)$ with $k_2 = \bQ[x]/(x^3-x^2-2x+1)$, and $\ell$ one of the fields from the above list, by $\cE_j$, $4\leqslant j\leqslant 10$.   
Note that three pairs belonging to the above two lists coincide with pairs of number fields in [PY1], \S 8:
$\cE_4=\cC_{31}$, $\cE_5=\cC_{32}$ and $\cE_1=\cC_{33}.$ 

\vskip1mm

Using the value of $\zeta_{\ell |k}(-2)$ given in the last column of the above two tables and the values of $\zeta_k(-1)$ and $\zeta_k(-3)$ given below for $k =k_1$ and $k_2$, we compute $\cR = 2^{-9}\zeta_k(-1))\zeta_{\ell |k}(-2)\zeta_k(-3)$ for each of the ten pairs $\cE_j$, $j\leqslant 10$. We find that the numerator of none of them is a power of $2$. Proposition 1 then implies that $d$ cannot be $3$ either.$$\zeta_{k_1}(-1)=-1/9,\ \  \zeta_{k_1}(-3)=199/90; \ \ \ \zeta_{k_2}(-1)=-1/21,\ \  \zeta_{k_2}(-3)=79/210.$$  

\vskip1mm

The following is a summary of what we have proved above.
\begin{prop} (i) If $n\geqslant 5$, then $d =1$, i.e., $k =\bQ$.\\
(ii) If $n=3$, then $d\leqslant 2.$ 
\end{prop}

\vskip3mm
\begin{center}
{\bf 9. $G$ of type $^2A_n$ with $n>1$ odd and $k= \bQ$}  
\end{center}
\vskip4mm

\ni{\bf 9.1.} 
 We shall assume in the sequel that $k = \bQ$ which according to Proposition 3 is the case if $n\geqslant 5$. 
 Then $r =1$ and   $\ell=\bQ(\sqrt{-a})$ for some square-free positive integer $a.$
By setting $d = 1$ and $D_k =1$ in bounds (61) and (62) we obtain

\begin{eqnarray*}
D_\ell &\leqslant& \kappa_1(n,h_{\ell,n+1}):=\ll{\mathfrak p}_1(n,1,1,h_{\ell,n+1})\lr.\\ 
D_{\ell}&\leqslant& \kappa_2(n,R_{\ell}/w_{\ell},\delta):= \ll{\mathfrak p}_2(n,1,1,R_{\ell}/w_{\ell},\delta)\lr. \\
\end{eqnarray*}
\ni{\bf 9.2.}
We easily see that for fixed $\delta \,(\geqslant 0.02)$ and $n$, $\kappa_2 $ decreases as 
$R_{\ell}/w_{\ell}$  increases, and for fixed $\delta\geqslant \,(0.02)$ and $R_{\ell}/w_{\ell}$, $\kappa_2$ decreases as $n$ increases  provided $n\geqslant 7$. 
Since the regulator of a complex quadratic field is $1,$ and $w_{\ell}=2$ for any complex quadratic field $\ell$ different from $\bQ(\sqrt{-3})$ and $\bQ(\sqrt{-1})$, $R_{\ell}/w_{\ell} = 1/2$ for all complex quadratic $\ell$ with $D_{\ell}>4$. 
Now for $n\geqslant 17$, as $\kappa_2(n, 1/8, 1.8)\leqslant \kappa_2 (17,1/8,1.8)=2,$ and there is no complex quadratic number field with discriminant $\leqslant 2$, {\it we conclude that $n\leqslant 15$}.  For $n=15,$ unless $\ell=\bQ(\sqrt{-3})$,
we know that $D_\ell\leqslant \kappa_2 (15,1/4,1.6)=3.$  Hence, if $n = 15$, $\ell=\bQ(\sqrt{-3})$.
For odd integers $n$ between $3$ and $13,$ unless $\ell=\bQ(\sqrt{-3})$ or $\bQ(\sqrt{-1}),$ with $D_\ell=3$ and $4$ respectively,
we can use the bound $D_\ell\leqslant\kappa_2 (n,1/2,\delta)$, with $\delta$ as indicated below, to obtain:
$$\begin{array}{cccccccc}
n:&&13&11&9&7&5&3\\
\delta:&&1.3&1&0.9&0.7&0.5&0.26\\
D_\ell\leqslant&&4&6&10&21&68&2874. 
\end{array}$$
\vskip1mm

\ni{\bf 9.3.} We will now improve the bound for the discriminant $D_\ell$ in case $n =3$. From the Table t20.001-t20.002 of [1] we see that the class number of every complex
quadratic number field $\ell$ with $D_\ell\leqslant 2874$ is  
$\leqslant  76$, and hence,   
$h_{\ell,4}\leqslant64$.  So we obtain the bound $D_\ell\leqslant \kappa_1(3,64)=1926.$
We can improve this bound further as follows.  From Table t20.001 we see that $h_\ell\leqslant 52,$ and hence $h_{\ell,4}\leqslant 32$, for all complex quadratic $\ell$ with $D_{\ell}\leqslant 1926$. Now we observe that, for $n =3$ and $k =\bQ$, combining equations (3), (4), and using the bounds (54), $(n+1)^r\mu(\cG/\Gamma)\leqslant 1$, $e'(P_v)\geqslant 1$ for all $v\in V_f$, and $e'(P_v)> n+1$ for all $v\in \cT$ (2.10), we get the following upper bound for $D_{\ell}$:
\begin{eqnarray}
D_{\ell}&<&\ll\Big[h_{\ell,4}\cdot \{ 
2\prod_{j=1}^{3}\frac{(2\pi)^{j+1}}{j!}\}\cdot\frac1{\zeta_{\bQ}(2)\zeta_{\ell |\bQ}(3)\zeta_{\bQ}(4)}
\Big]^{\frac{2}{5}}\lr\\
&\leqslant& \ll\Big[h_{\ell,4}\cdot \{ 
2\prod_{j=1}^{3}\frac{(2\pi)^{j+1}}{j!}\}\cdot\frac1{\zeta_\bQ(2)^{1/2}\zeta_\bQ(4)}
\Big]^{\frac{2}{5}}\lr=:\widetilde\kappa_1(h_{\ell,4})
\end{eqnarray}
where we have used the fact that $\zeta_{\bQ}(2)^{1/2}\zeta_{\ell|\bQ}(3)>1$ (see Lemma 1 in [PY2]).   Hence we conclude that
$D_\ell\leqslant \widetilde\kappa_1(32)\leqslant 1363.$

\ms
\ni{\bf 9.4.}  We can improve the bounds for $D_{\ell}$ for  $15>  n\geqslant 5$ as follows.  
 From the table of complex
quadratics in [C], we know that
$h_\ell\leqslant 5$  for $D_\ell\leqslant 68.$  Hence
$h_{\ell,n+1}\leqslant 5$ for $n\geqslant5.$

We now compute the values of $\kappa_1(n,j)$ for $5\leqslant n< 15$ and $1\leqslant j\leqslant 5.$

$$\begin{array}{cccccccc}
&\kappa_1(n,1)&\kappa_1(n,2)&\kappa_1(n,3)&\kappa_1(n,4)&\kappa_1(n,5)\\
n=5&47&52&55&57&59\\
n=7&18&19&20&20&20\\
n=9&10&10&10&10&10\\
n=11&6&6&6&6&6\\
n=13&4&4&4&4&4.
\end{array}$$

Comparing the above table with the table of complex quadratic number fields (cf. [C]) in terms of discriminants and class number,
we obtain the following possibilities for $D_{\ell}$ and $a$ (recall that $\ell = \bQ(\sqrt{-a})$):

$$\begin{array}{ccc}
n&D_{\ell}&a\\
15&3&3\\
13&3,4&3,1\\
11&3,4&3,1\\
9&3,4,7,8&3,1,7,2\\
7&3,4,7,8,11,15&3,1,7,2,11,15\\
5&3,4,7,8,11,15,19,20,23&3,1,7,2,11,15,19,5,23\\
&24,31,35,39,40,43&
6,31,35,39,10,43\\
&47,51,52,55,56&47,51,13,55,14.
\end{array}$$
In the above we have used the fact
 that for $\ell = \bQ (\sqrt{-a} )$ , where $a$ is a 
square-free positive integer, $D_{\ell} =a$ if $a\equiv 3$ (${\rm mod}\ 4$), and $D_{\ell} =4a$ otherwise.

\vskip2mm
\ni{\bf 9.5.}  To prove Theorem 2 (stated in the Introduction) we  compute $\cR$ in each of the cases occurring  in the second table of 9.4 using the following values of $\zeta:=\zeta_{\bQ}$ and  $\zeta_{\ell |\bQ}$. 

$$
\begin{array}{ccccccccl}
j:&-1&-3&-5&-7&-9&-11&-13&-15\\
\zeta(j):&
-1/12&1/120&-1/252&1/240&-1/132&691/32760&-1/12&3617/8160.\end{array}$$

\ni Listed below are the values of $\zeta_{\ell|k}$, for $(k,\ell)=(\bQ,\bQ(\sqrt{-3}))$,

$$\begin{array}{ccccccc}
\zeta_{\ell|\bQ}(-2)&\zeta_{\ell|\bQ}(-4)&\zeta_{\ell|\bQ}(-6)&\zeta_{\ell|\bQ}(-8)&\zeta_{\ell|\bQ}(-10)&\zeta_{\ell|\bQ}(-12)&\zeta_{\ell|\bQ}(-14)\\
-2/9&2/3&-14/3&1618/27&-3694/3&111202/3&13842922/9,
\end{array}
$$ 
and the values of $\zeta_{\bQ(\sqrt{-a})|\bQ }$ required to compute $\cR$ for $5\leqslant n\leqslant 13$, $a\ne 3$, are given below:
$$\begin{array}{ccccccccl}
a&\zeta_{\ell|\bQ}(-2)&\zeta_{\ell|\bQ}(-4)&\zeta_{\ell|\bQ}(-6)&\zeta_{\ell|\bQ}(-8)&\zeta_{\ell|\bQ}(-10)&\zeta_{\ell|\bQ}(-12)\\
1&-1/2&5/2&-61/2&-1385/2&-50521/2&2702765/2\\
7&-16/7&32&-1168&565184/7&&\\
2&-3&57&-2763& 250737&&\\
11&-6&2550/11&-21726&&&\\
15& -16&992& -165616.&&&
\end{array}
$$ 

$$\begin{array}{ccccccccccccccccc}
a:&19&5&23&6&31&35&39\\
\zeta_{\ell |\bQ}(-2):&-22&-30&-48&-46&-96&-108&-176\\
\zeta_{\ell |\bQ}(-4):&2690&3522&6816&7970&25920&42372&73120.\\
\end{array}$$
$$\begin{array}{ccccccccccccccl}
a:&10&43&47&51&13&55&14\\
\zeta_{\ell |\bQ}(-2):&-158&-166&-288&-268&-302&-400&-396\\
\zeta_{\ell |\bQ}(-4):&79042&106082&169920&229700&257314&341984&362340.\\
\end{array}$$
\vskip1mm

Explicit computation of $\cR$ in each of the above cases shows that for every odd integer $n>7$, the numerator of $\cR$ has a prime divisor which does not divide $n+1$. In view of Propositions 1 and 3, this proves Theorem 2 for $n>7$. 

\vskip1mm

For $n = 7, \,5$, we list below the value of $\cR$ for those $a$  in the second table in 9.4 for which the prime divisors of the numerator of $\cR$ divide $n+1$.
 $$\begin{array}{cccccl}
n&a&\zeta_{\ell|\bQ}(-2)&\zeta_{\ell|\bQ}(-4)&\zeta_{\ell|\bQ}(-6)&\cR\\
7&3&-2/9& 2/3&-14/3&1/16124313600=1/(2^{15}\cdot3^9\cdot5^2)\\
5&3&-2/9& 2/3&&1/78382080=1/(2^{10}\cdot3^7\cdot5\cdot7)\\
5&1&-1/2&5/2&&1/9289728=1/(2^{14}\cdot3^4\cdot7)\\
5&7&-16/7&32&&1/158760=1/(2^3\cdot 3^4\cdot 5\cdot 7^2)\\
5&31&-96&25920&&3/14.\\
\end{array}
$$ 
We need to consider only the $a$ appearing in the above table.  
\vskip2mm

\ni{\bf 9.6.} In our treatment of groups of type $^2A_n$, with $n$ odd, we have not so far made use of the assumption that $\Gamma$ is cocompact, or, equivalently, $G$ is anisotropic over $k$, see 1.5. We will now use the fact that  $G$ is anisotropic over $k=\bQ$ to exclude $n = 7,\, 5$. This will complete our proof of Theorem 2. 
\vskip1mm

From the well-known description of absolutely simple simply connected $\bQ$-groups of type $^2A_n$ we know that there is a division algebra $\cD$ with center $\ell$ and of degree $\mathfrak{d} = \sqrt{[\cD:\ell]}$, $\mathfrak{d}|(n+1)$, $\cD$ given with an involution $\sigma$ of the second kind, and a nondegenerate hermitian form $h$ on $\cD^{(n+1)/\mathfrak{d}}$ defined in terms of the involution $\sigma$, so that $G$ is the special unitary group ${\mathrm{SU}}(h)$ of $h$.
\vskip 1mm

If $\cD =\ell$, then $h$ is an hermitian form on $\ell^{n+1}$ such that the quadratic form $q$ on the $2(n+1)$-dimensional $\bQ$-vector space $V = \ell^{n+1}$ defined by $$q(v) = h(v,v) \ \ \ {\mathrm{for}}\ \ v\in V,$$ is isotropic over $\bR$ (since $G$ is isotropic over $\bR$, i.e., $G(\bR)$ is noncompact). Then as $n\geqslant 3$, $q$ is isotropic over $\bQ$ by Meyer's theorem and hence G is isotropic over $\bQ$. But this is not the case. Therefore, $\cD \ne \ell$, i.e., $\cD$ is a noncommutative division algebra of degree $\mathfrak{d}>1$.  
\vskip1mm

Using the structure of the Brauer group of a global field, we see that there exists at least one prime $p$ which splits over $\ell$ such that $\bQ_p\otimes_{\bQ}\cD= (\bQ_p\otimes_{\bQ}\ell)\otimes_{\ell}\cD$ is isomorphic to $M_m(\mathfrak{D}_p)\times M_m(\mathfrak{D}_p^o)$, where $\mathfrak{D}_p$ is a  noncommutative central division algebra over $\bQ_p$ of degree $\mathfrak{d}_p>1$, $\mathfrak{D}_p^o$ is its opposite, and $m = (n+1)/\mathfrak{d}_p$. The involution $\sigma$ interchanges the two factors of $M_m(\mathfrak{D}_p)\times M_m(\mathfrak{D}_p^o)$, and hence, $G(\bQ_p) \cong {\mathrm{SL}}_m(\mathfrak{D}_p)$. 
\vskip1mm

In the rest of this section $n$ is either $7$ or $5$, $a$ and $\cR$ are as in the last table of 9.5. 
The nonarchimedean place of $\bQ$ corresponding to a prime $p$ will be denoted by $p$.  Now let $p$ be a prime which splits in $\ell = \bQ(\sqrt{-a})$ and $\mathfrak{D}_p$ is  a noncommutative division algebra with center $\bQ_p$. 
Then (see the computation in 2.3(ii) of  [PY2]) $e'(P_p)$ is an integral multiple of $f_7(p):=(p-1)(p^3-1)(p^5-1)(p^7-1)$ if $n =7$, and it is an integral multiple of either $f_5(p):=(p-1)(p^3-1)(p^5-1)$ or $g_5(p):=(p-1)(p^2-1)(p^4-1)(p^5-1)$ if $n =5$.  

\vskip1mm

Now let $\cT$ be as in 2.10. Recall that for every prime $q$, $e'(P_q)$ is an integer, and for $q\in \cT$, $e'(P_q)>n+1$. Also recall that $\mu(\cG/\Gamma)$ is a submultiple of $1/(n+1)$, and 
\begin{equation}
\mu(\cG/\Gamma)= \frac{\cR\prod e'(P_q)}{[\Gamma:\Lambda]}=\frac{\cR e'(P_p)\prod_{q\ne p}e'(P_q)}{[\Gamma:\Lambda]}.
\end{equation} 
As every prime divisor of $[\Gamma:\Lambda]$ divides $n+1$, we conclude that every prime divisor of the numerator of $\cR e'(P_p)$ divides $n+1$.   Also since $[\Gamma:\Lambda]\leqslant 2 h_{\ell, n+1}(n+1)^{1+\#\cT}$ (cf.\,(54)), we see that \begin{equation}\frac{\cR e'(P_p)}{2h_{\ell, n+1}(n+1)^2}\leqslant \mu(\cG/\Gamma)\leqslant \frac1{(n+1)},\end{equation} and hence, \begin{equation}\cR e'(P_p)\leqslant 2h_{\ell,n+1}(n+1).\end{equation} Now we note that the class number of the complex quadratic field $\ell = \bQ(\sqrt{-a})$, for $a =3,\, 1,\, 7$ is $1$, and for $a=31$ the class number is $3$.  The first two primes $\{p_1,\,p_2\}$ which split in $\bQ(\sqrt{-a})$ are $\{7,\,13\}$, $\{5,\,13\}$, $\{2,\,11\}$ and $\{2, \, 5\}$ for $a = 3, \,1,\,7$ and $31$ respectively.  Let $\cR$ be as in the last column of the last table of 9.5. By direct computations we see that if $n =7$ and $a =3$, $\cR e'(P_p)\geqslant \cR f_7(7)>  16$, and if $n =5$ and $a =31$,  both $\cR f_5(2)$ and $\cR g_5(2)$ are larger than $36$. On the other hand, if $n =5$ and $a = 3,\,1$ or $7$, both $\cR f_5(p_2)$ and $\cR g_5(p_2)$ are larger than $12$, and at least one prime divisor of  the numerator of $\cR f_5(p_1)$ and $\cR g_5(p_1)$ is different from $2$ and $3$. We conclude from these observations that $n$ cannot be $5$ or $7$. Thus we have proved Theorem 2.  
\vskip3mm

\ni{\bf Corrections in} [PY2]: (i)  In line 11 on page 381 and in the last two lines on page 386,  ``$\chi(\Gamma)$''  and ``$\chi(\Lambda)$''should be replaced with ``$\vert\chi(\Gamma)\vert$'' and ``$\vert\chi(\Lambda)\vert$'' respectively. (ii) In the statement of Theorem 2 on page 402, ``$\chi(X_u)/n$" should be replaced with ``$\chi(X_u)$".
\vskip1mm

We note that a revised version of [PY1] which incorporates corrections and additions given in the ``Addendum'' has recently been posted on the arXiv. 

\vskip3mm    
\ni{\bf Acknowledgments.} We thank Gunter Malle for providing us lists of number fields used in this paper.
\vskip1mm

The first-named author was supported by the Humboldt Foundation and the 
NSF (grant DMS-1001748).  
The second-named author received partial support from the NSA.  The paper was completed while the
second-named author visited the Institute of Mathematics of the University of Hong Kong, to which he
would like to express his gratitude. 

\bs
\centerline{\bf References}

\vskip4mm

\ni[B] A.\,Borel, {\it Stable real cohomology of arithmetic groups}, 
Ann.\,Sci.\,Ec.\,Norm.\,Sup.\,(4) {\bf 7}(1974), 235-272. 
\vskip1.5mm

\ni[BP] A.\,Borel and G.\,Prasad, {\it Finiteness theorems for discrete subgroups of bounded covolume in semisimple groups.} Publ.\,Math.\,IHES No.\,{\bf 69}(1989), 119--171.

\vskip1.5mm
\ni [C]  H.\,Cohen, {\it A course in computational algebraic number theory.} Graduate Texts in Mathematics, 138. Springer-Verlag, Berlin, 1993.

\vskip1.5mm

\ni[CS] D.\,I.\,Cartwright and Tim Steger, {\it Enumeration of the 50 fake projective planes.} C.\,R.\,Acad. Sc.\,Paris, Ser.\,I\, {\bf 348}(2010), 11-13.
\vskip1.5mm

\ni[F] E. Friedman, {\it Analytic formulas for the regulator of a number field.}
Invent.\,Math., {\bf 98}(1989), 599--622.

\vskip1.5mm

\ni[Gi] P.\,Gille, {\it The Borel-de Siebenthal's theorem}, preprint (2008).
\vskip1.5mm

\ni[GM] A.\,S.\,Golsefidy and A.\,Mohammadi, {\it Discrete subgroups acting transitively on vertices of a Bruhat-Tits building}, Duke Math.\,J.\,{\bf 161}(2012), 483-544.
\vskip1.5mm
 
\ni [H] S.\,Helgason, {\it Differential geometry, Lie groups, and symmetric spaces}, Academic press, 1978.

\vskip1.5mm

\ni[KK] V.\,Kharlamov and V.\,Kulikov, {\it On real structres on rigid surfaces.}  Izv.\,Math.\,{\bf 66}(2002), 133-150.
\vskip1.5mm

\ni[Marg] G.A.\,Margulis, {\it Discrete subgroups of semi-simple Lie groups}, Springer-Verlag, Heidelberg (1991).

\vskip1.5mm
\ni [Mart] Martinet, J.,  {\it Petits discriminants des corps de nombres.} Number theory days (Exeter, 1980), 151--193, London Math.\:Soc.\:Lecture Note Ser., 56, Cambridge Univ.\:Press, Cambridge-New York, 1982.

\vskip1.5mm
\ni [Mu] D.\,Mumford, {\it An algebraic surface with $K$ ample}, $K^2=9$, $p_g=q=0.$ Amer.\,J.\,Math. {\bf 101}(1979), 233--244.

\vskip1.5mm

\ni[N] W.\,Narkiewicz, {\it Elementary and analytic theory of algebraic numbers,} third edition. Springer-Verlag, New York (2000).

\vskip1.5mm

\ni[O] A.\,M.\,Odlyzko, {\it Discriminant bounds,} unpublished, available
from:

 http://www.dtc.umn.edu/$\sim$odlyzko/unpublished/index.html.
\vskip1.5mm

\ni[P]  G.\,Prasad, {\it Volumes of $S$-arithmetic quotients of semi-simple groups.} Publ.\,Math. IHES No.\,{\bf 69}(1989), 91--117.
\vskip1.5mm

\ni[PR] G.\,Prasad and A.S.\,Rapinchuk, {\it Weakly commensurable arithmetic groups and isospectral locally symmetric spaces}, Publ.\,Math.\,IHES, {\bf 119}(2009), 113-184. .

\vskip1.5mm
\ni[PY1] G.\,Prasad and S-K.\,Yeung, {\it Fake projective planes}. Inv.\,Math.\,{\bf 168}(2007), 321-370. {\it Addendum}, ibid {\bf 182}(2010), 213-227.

\vskip1.5mm
\ni
[PY2] G.\,Prasad and S-K.\,Yeung, {\it Arithmetic fake projective spaces and arithmetic fake Grassmannians}, American J.\,Math.\,{\bf 131}(2009), 379-407.

\vskip1.5mm
\ni
[Sel] S. Selmane, {\it  Quadratic extensions of totally real quintic fields.}  Math.\,Comp.\,{\bf 70} (2001), 837-843.

\vskip1.5mm

\ni[Ser] J-P.\:Serre, {\it Cohomologie des groupes discrets,} in Annals of 
Math.\,Studies {\bf 70}. Princeton U.\,Press, Princeton (1971).
\vskip1.5mm

\ni[Ta] M.\,Takeuchi, {\it On the fundamental group and the group of isometries of a symmetric space}, J.\,Fac.\,Sci.\,Univ.\,Tokyo {\bf 10}(1964), 88-123.

\vskip1.5mm

\ni[Ti1] J.\,Tits, {\it Classification of algebraic semisimple groups.} Algebraic Groups and Discontinuous
Subgroups. Proc.\,A.M.S.\,Symp.\,Pure Math.\:{\bf 9}(1966) pp. 33--62. 
\vskip1.5mm

\ni[Ti2] J.\,Tits, {\it Reductive groups over local fields.} Proc.\,A.M.S.\,Symp.\,Pure 
Math. {\bf 33}(1979), Part I, 29--69.
 \vskip1.5mm

\ni[W] L.\,C.\,Washington, {\it Introduction to Cyclotomic Fields}, 2nd edition, Grad.\,Texts Math., vol.\,83, Springer, New York (1997).
\vskip1.5mm
 
\ni[Z] R.\,Zimmert, {\it Ideale kleiner Norm in Idealklassen und eine Regulatorabsch\"atzung}, Inv.\,Math. {\bf {62}}(1981), 
367-380.
\vskip1.5mm 

\ni[1] The Bordeaux Database, Tables obtainable from:

    ftp://megrez.math.u-bordeaux.fr/pub/numberfields/.
\vskip6mm

\vskip4mm\end{document}